\documentclass[12pt,leqno]{amsart}

\usepackage[colorlinks=true, pdfstartview=FitV, linkcolor=blue,%
citecolor=blue, urlcolor=blue]{hyperref}

\usepackage{amssymb,amsmath,amsthm}
\usepackage[all]{xy}

\usepackage{enumerate}
\usepackage{mathrsfs}
\usepackage{graphics}
\usepackage{verbatim}

\usepackage[usenames,dvipsnames]{xcolor}

\usepackage{marginnote}
\usepackage{accents}
\usepackage{xspace}

\numberwithin{equation}{section}

\newenvironment{rouge}
{\relax\color{red}}
{\hspace*{.3ex}\relax}
\newcommand{\berm}{\begin{rouge}{}\marginnote{\fbox{\scshape\lowercase{M}}}{}}
\newenvironment{bleu}
{\relax\color{PineGreen}}
{\hspace*{.3ex}\relax}
\newcommand{\beb}{\begin{bleu}}
\newcommand{\eb}{\end{bleu}}
\newdir{ >}{{}*!/-10pt/@{>}}

\renewcommand{\emptyset}{\varnothing}

\newcommand{\nc}{\newcommand}
\nc{\on}{\operatorname}

\newlength{\my}
\nc{\noi}{\noindent}

\renewcommand{\Re}{\operatorname{Re}}

\newtheorem{theorem}{Theorem}[section]
\newtheorem{proposition}[theorem]{Proposition}
\newtheorem{lemma}[theorem]{Lemma}
\newtheorem{corollary}[theorem]{Corollary}

\theoremstyle{definition}

\newtheorem{definition}[theorem]{Definition}
\newtheorem{notation}[theorem]{Notation}
\newtheorem{example}[theorem]{Example}
\newtheorem{examples}[theorem]{Examples}
\newtheorem{remark}[theorem]{Remark}

\newtheorem{conjecture}[theorem]{Conjecture}

\nc{\Rem}{\begin{remark}}
\nc{\enrem}{\end{remark}}
\nc{\Conj}{\begin{conjecture}}
\nc{\enconj}{\end{conjecture}}
\nc{\Th}{\begin{theorem}}
\nc{\enth}{\end{theorem}}
\nc{\Lemma}{\begin{lemma}}
\nc{\enlemma}{\end{lemma}}
\nc{\Cor}{\begin{corollary}}
\nc{\encor}{\end{corollary}}
\nc{\Def}{\begin{definition}}
\nc{\edf}{\end{definition}}
\nc{\Ex}{\begin{example}}
\nc{\enex}{\end{example}}
\nc{\Exs}{\begin{examples}}
\nc{\enexs}{\end{examples}}

\renewcommand{\ge}{\geqslant}
\renewcommand{\le}{\leqslant}

\nc{\RR}{\mathrm{R}}
\nc{\LL}{\mathrm{L}}

\newcommand{\C}{{\mathbb{C}}}

\newcommand{\R}{{\mathbb{R}}}

\newcommand{\Z}{{\mathbb{Z}}}

\newcommand{\BBP}{{\mathbb P}}

\newcommand{\BBD}{{\mathbb D}}

\newcommand{\PR}{{\BBP^1(\R)}}

\newcommand{\cor}{\C}

\newcommand{\icor}{\cor^\sub}

\newcommand{\JD}{\D^\sub}

\newcommand{\iC}{\C^\sub}
\newcommand{\iD}{\D^\sub}

\nc{\tA}{\mathcal{A}}
\nc{\rA}{\mathscr{A}}

\def\D{\mathscr{D}}
\def\DA{\D^\tA}

\def\sha{\mathscr{A}}

\def\shd{\mathscr{D}}
\def\she{\mathscr{E}}
\def\shf{\mathscr{F}}

\def\shi{\mathscr{I}}

\def\shl{\mathscr{L}}
\def\shm{\mathscr{M}}
\def\shn{\mathscr{N}}
\def\sho{\mathscr{O}}

\def\shs{\mathcal{S}}

\renewcommand{\ker}{\operatorname{Ker}}

\newcommand{\pt}{\mathrm{pt}}

\newcommand{\into}{\hookrightarrow}

\def\monoto{{\rightarrowtail}}

\renewcommand{\to}[1][]{\xrightarrow[]{#1}}
\newcommand{\from}[1][]{\xleftarrow[]{#1}}
\newcommand{\isoto}[1][]{\xrightarrow[#1]%
{{\raisebox{-.6ex}[0ex][-.6ex]{$\mspace{1mu}\sim\mspace{2mu}$}}}}

\newcommand{\To}[1][\rule{1ex}{0pt}]{\xrightarrow{\hs{.6ex}#1\hs{.6ex}}}

\newcommand{\muHom}[1][]{\mathrm{Hom}^\mu_{\raise1.5ex\hbox to.1em{}#1}}
\newcommand{\Hom}[1][]{\mathrm{Hom}_{\raise1.5ex\hbox to.1em{}#1}}
\newcommand{\RHom}[1][]{\RR\mathrm{Hom}_{\raise1.5ex\hbox to.1em{}#1}}
\newcommand{\Ext}[2][]{\mathrm{Ext}_{\raise1.5ex\hbox to.1em{}#1}^{#2}}
\renewcommand{\hom}[1][]{{\mathscr{H}\mspace{-4mu}om}_{\raise1.5ex\hbox to.1em{}#1}}
\newcommand{\rhom}[1][]{{\RR\mathscr{H}\mspace{-3mu}om}_{\raise1.5ex\hbox to.1em{}#1}}
\newcommand{\rhomc}[1][]
{{\mathscr{H}\mspace{-3mu}om}^*_{\raise1.5ex\hbox to.1em{}#1}}

\newcommand{\cihom}[1][]
{{\mathscr{I}\mspace{-3mu}}{hom}^+_{\raise1.5ex\hbox to.1em{}#1}}

\nc{\ihom}[1][]{{\shi\mspace{-3mu}hom}_{\raise1.5ex\hbox to.1em{}#1}}
\nc{\rihom}[1][]{{\mspace{2mu}\mathrm{R}\shi\mspace{-3mu}hom}_{\raise1.5ex\hbox to.1em{}#1}}
\nc{\fihom}[1][]{{\shi\mspace{-3mu}hom}^{\mathsf{E}}_{\raise1.5ex\hbox to.1em{}#1}}
\nc{\FHom}[1][]{{\mathrm{RHom}^{\mathrm{E}}_{\raise1.5ex\hbox to.1em{}#1}}}
\nc{\fhom}[1][]{{\mathscr{H}%
\mspace{-3mu}om}^{\mathsf{E}}_{\raise1.5ex\hbox to.1em{}#1}}
\nc{\Endom}[1][]{{\she\mspace{-3mu}nd}_{\raise1.5ex\hbox to.1em{}#1}}

\newcommand{\eeim}[1]{{#1}_{!!}\ms{3mu}}
\newcommand{\ssubset}{\subset\mspace{-3mu}\subset}

\nc{\Tam}{{\mathrm{E}}}

\newcommand{\ext}[2][]{{\mathscr{E}xt}_{\raise1.5ex\hbox to.1em{}#1}^{#2}}
\newcommand{\Tor}[2][]{\mathrm{Tor}^{\raise1.5ex\hbox to.1em{}#1}_{#2}}
\newcommand{\tens}[1][]{\mathbin{\otimes_{\raise1.5ex\hbox to-.1em{}{#1}}}}
\newcommand{\ltens}[1][]{\mathbin{\overset{\mathrm{L}}\tens}_{#1}}

\newcommand{\etens}{\mathbin{\boxtimes}}

\newcommand{\Endo}[1][]{\mathrm{End}_{\raise1.5ex\hbox to.1em{}#1}}

\newcommand{\Aut}[1][]{\mathrm{Aut}_{\raise1.5ex\hbox to.1em{}#1}}
\newcommand{\sect}{\Gamma}

\newcommand{\ctens}[1][]{\mathbin{\overset{+}\tens}_{#1}}

\newcommand{\VV}{{\mathsf{V}}}
\newcommand{\VVd}{{\mathsf{V}^*}}

\newcommand{\WW}{{\mathbb{V}}}

\newcommand{\oim}[1]{{#1}_*}
\newcommand{\eim}[1]{{#1}_{\ms{1mu}!}\ms{3mu}}
\newcommand{\roim}[1]{\RR{#1}_*}
\newcommand{\reim}[1]{\RR{#1}_{\ms{1mu}!}\ms{3mu}}

\newcommand{\reeim}[1]{\RR{#1}_{\mspace{1mu}!!}\ms{3mu}}
\newcommand{\opb}[1]{#1^{-1}}

\newcommand{\epb}[1]{#1^{\ms{3mu}!}\ms{3mu}}

\newcommand{\Dtens}[1][]{\overset{\mathrm{D}}\otimes_{\raise1.5ex\hbox to-.1em{}#1}}
\newcommand{\Detens}[1][]{\overset{\mathrm{D}}\etens_{\raise1.5ex\hbox to-.1em{}#1}}
\newcommand{\Ddual}{{\BBD}}

\newcommand{\Deim}[1]{{\mathrm{D}}{\mspace{-1mu}#1_{\mspace{1mu}!}}}
\newcommand{\Dopb}[1]{{\mathrm{D}}{#1}^{*}}
\newcommand{\Doim}[1]{{\mathrm{D}}{#1}_{*}}

\newcommand{\good}{\mathrm{good}}
\newcommand{\qgood}{\mathrm{q\text-good}}

\newcommand{\Tmp}{\mathsf{T}}
\nc{\rE}{\mathrm{E}}
\nc{\enh}{\Tmp}

\newcommand{\Toim}[1]{\RR\ms{.7mu}{{#1}_\R}_{\ms{.5mu}*}}

\newcommand{\Topb}[1]{({#1}_\R)^{-1}}
\newcommand{\Tepb}[1]{({#1}_\R)^{\,!}}

\nc{\EF}[1][]{{}^{\enh}\mspace{-3mu}\shf_{#1}}
\nc{\EFa}[1][]{{}^{\enh}{\mspace{-3mu}\shf^a_{#1}}}
\nc{\FS}[1][]{{}^{\mathrm{S}}{\mspace{-3mu}\shf_{#1}}}
\nc{\FSa}[1][]{{}^{\mathrm{S}}{\mspace{-3mu}\shf^a_{#1}}}
\nc{\Leg}[1][]{{\mathrm{Conv}}{(#1)}}
\nc{\dom}{\mathrm{dom}}
\nc{\domo}{\dom^\circ}

\nc{\hol}{\mathrm{hol}}

\nc{\Mod}{\mathrm{Mod}}
\nc{\rh}{\mathrm{rh}}
\newcommand{\md[1]}{\Mod({#1})}

\newcommand{\mdc[1]}{{\Mod_{\ms{3mu}\coh}}({#1})}

\newcommand{\mdhol[1]}{\Mod_{\ms{3mu}\hol}({#1})}

\newcommand{\mdrc[1]}{\Mod_{\Rc}({#1})}
\newcommand{\mdrcc[1]}{\Mod^{\mathrm{c}}_{\Rc}({#1})}

\nc{\sHH}{\mathscr{H}\mspace{-4mu}\mathscr{H}}

\nc{\sMH}{\mathscr{M}\mspace{-4mu}\mathscr{H}}

\newcommand{\eqdot}{\mathbin{:=}}
\newcommand{\seteq}{\mathbin{:=}}

\newcommand{\cl}{\colon}
\newcommand{\scbul}{{\,\raise.4ex\hbox{$\scriptscriptstyle\bullet$}\,}}

\newcommand{\tw}[1]{\widetilde{#1}}

\newcommand{\twX}{{\widetilde{X}}}
\nc{\twL}{\tw{L}}

\newcommand{\ol}{\overline}
\newcommand{\bl}{\bigl(}
\newcommand{\br}{\bigr)}
\newcommand{\lp}{{\rm(}}
\newcommand{\rp}{{\rm)}}

\newcommand{\Cc}{{\C\text{-c}}}
\newcommand{\Rc}{{\ms{1mu}\R\text{-c}}}

\newcommand{\Sol}[1][X]{\mspace{1mu}{\shs\mspace{-2.5mu}\mathit{ol}}_{#1}}

\nc{\be}{\begin{enumerate}}
\nc{\ee}{\end{enumerate}}
\newcommand{\bnum}{\begin{enumerate}[{\rm(i)}]}
\newcommand{\enum}{\end{enumerate}}
\newcommand{\banum}{\begin{enumerate}[{\rm(a)}]}
\newcommand{\eanum}{\end{enumerate}}

\newcommand{\bna}{\be[{\rm(a)}]}

\newenvironment{myequation}
{\relax\setlength{\arraycolsep}{1pt}\begin{eqnarray}}
{\end{eqnarray}}
\newenvironment{myequationn}
{\relax\setlength{\arraycolsep}{1pt}\begin{eqnarray*}}
{\end{eqnarray*}}

\nc{\eq}{\begin{myequation}}
\nc{\eneq}{\end{myequation}}
\nc{\eqn}{\begin{myequationn}}
\nc{\eneqn}{\end{myequationn}}

\newenvironment{myarray}[2][c]{\relax\setlength{\arraycolsep}{1pt}

\begin{array}[#1]{#2}}{\end{array}\relax}

\newcommand{\ba}{\begin{myarray}}
\newcommand{\ea}{\end{myarray}}

\newcommand{\set}[2]{\left\{#1 \mathbin{;} #2 \right\}}

\nc{\Proof}{\begin{proof}}
\nc{\QED}{\end{proof}}
\nc{\Prop}{\begin{proposition}}
\nc{\enprop}{\end{proposition}}

\nc{\rop}{{\ms{3mu}\mathrm{op}}}
\nc{\op}{\rop}
\nc{\tot}{\mathrm{tot}}
\nc{\Op}[1][M]{{\mathrm{Op}_{#1}}}
\nc{\dist}{{\mathrm{dist}}}
\nc{\LocSyst}{{\mathrm{LocSyst}}}
\nc{\eu}{\mathrm{eu}}
\nc{\hh}{\mathrm{hh}}
\nc{\mueu}{{\mu\eu}}

\DeclareMathOperator{\id}{id}

\DeclareMathOperator{\chv}{char}

\newcommand{\Supp}{\on{Supp}}
\newcommand{\Der}[1][]{\mathsf{D}^{#1}}
\newcommand{\Derb}{\Der[\mathrm{b}]}

\newcommand{\SSupp}{\on{Sing\ms{.5mu}Supp}}

\newcommand{\Derp}{\Der[+]}
\newcommand{\Derm}{\Der[-]}

\newcommand{\coh}{\mathrm{coh}}

\newcommand{\dT}{{\dot{T}}}

\newcommand{\BDC}{\Derb}
\newcommand{\TDC}{\Derb}

\newcommand{\mop}{\mathrm{r}}

\newcommand{\OEn}{\sho^{\mspace{2mu}\enh}}
\newcommand{\DbT}{\Db^\Tmp}

\newcommand{\OvE}{\Omega^\enh}
\newcommand{\drE}{\mathcal{DR}^\Tmp}
\newcommand{\solE}{\mspace{1mu}\mathcal{S}ol^{\mspace{2mu}\Tmp}}

\nc{\wc}[1]{\overset{\mbox{$\scriptscriptstyle\vee$}}{#1}}
\nc{\field}{\cor}

\nc{\bVd}{{\ol \VVd}}
\nc{\bWd}{{{\ol \WW}^*}}

\newcommand{\fR}{{\R_\infty}}

\nc{\oX}{{\ol X}}
\nc{\oS}{\ol S}
\nc{\oY}{\ol Y}
\nc{\oL}{\ol L}
\nc{\oR}{{\ol\R}}
\nc{\Tl}{\mathrm{L^E}}
\nc{\Tr}{\mathrm{R^E}}

\nc{\sa}{\mathrm{sa}}

\newcommand{\Db}{{\mathcal D} b}
\newcommand{\Dbt}{{\mathcal D} b^{\mspace{2mu}\mathrm t}}

\newcommand{\Ot}[1][X]{\sho^{\mspace{2.5mu}{\mathrm t}}_{#1}}
\newcommand{\OO}[1][X]{\sho_{#1}}

\nc{\BRC}{{\R}\hbox{-}{\mathrm{Cons}}}
\nc{\Brc}{{\R}\hbox{-}{\mathrm{C}}}
\newcommand{\Ovt}{\Omega^{\mspace{1.5mu}{\mathrm t}}}

\newcommand{\dr}{\mathcal{DR}}
\newcommand{\drt}{\mathcal{DR}^{\mathrm t}}

\newcommand{\sol}{\mathcal Sol}

\newcommand{\solt}{\mathcal Sol^{\mspace{2.5mu}\mathrm t}}
\newcommand{\reghol}{\rh}

\newcommand{\indlim}[1][]{\mathop{\varinjlim}\limits_{#1}}
\newcommand{\sindlim}[1][]{\smash{\mathop{\varinjlim}\limits_{#1}}\,}

\newcommand{\prolim}[1][]{\mathop{\varprojlim}\limits_{#1}}

\newcommand{\inddlim}[1][]{\mathop{\text{\rm``{$\varinjlim$}''}}\limits_{#1}}
\newcommand{\sinddlim}[1][]{\smash{\mathop{\text{\rm``{$\varinjlim$}''}}\limits_{#1}}\,}

\newcommand{\dsum}[1][]{\mathbin{\oplus_{#1}}}
\nc{\eps}{\varepsilon}
\nc{\hs}{\hspace*}
\nc{\nn}{\nonumber}
\nc{\tM}{\widetilde{M}}
\nc{\h}{\mathbf{h}}
\nc{\tf}{\tilde{f}}
\nc{\trf}{{{}^{\mathrm{t}}\mspace{-3mu}f}}
\nc{\codim}{\on{codim}}
\nc{\lh}{\mathscr{H}}
\nc{\bwr}{\scalebox{1.1}{$\wr$}}
\nc{\dTi}{\dT^{*,\mathrm{in}}}
\nc{\Cd}{\mathrm{C}}
\nc{\tK}{\widetilde{K}}
\nc{\aMM}{a_{M\times M}}
\nc{\e}{\mspace{1mu}\mathrm{e}\mspace{1mu}}
\nc{\lan}{\langle}
\nc{\ran}{\rangle}
\nc{\la}{\lambda}

\newcommand{\Union}{\bigcup\limits}

\nc{\vphi}{\varphi}
\nc{\vep}{\varepsilon}

\nc{\At}{\tA_\twX}
\nc{\bu}{\boldsymbol{u}}
\nc{\bv}{\boldsymbol{v}}
\nc{\hbu}{\widehat{\boldsymbol{u}}}
\nc{\ex}{\mathrm{e}}
\nc{\vpi}{\varpi}
\newcommand{\soplus}{\mathop{\mbox{\small $\bigoplus$}}\limits}
\nc{\one}{\mathbf{1}}
\nc{\setp}[1]{\{#1\}}
\nc{\GL}{\mathrm{GL}}
\nc{\cI}{\mathrm{I}}

\nc{\vs}{\vspace*}

\nc{\wb}[1]{\mbox{$\rule[-1.1ex]{0ex}{2ex}#1$}}
\nc{\wwb}[1]{\mbox{$\rule[-1.8ex]{0ex}{3ex}#1$}}
\nc{\bpi}{\ol{\pi}}
\nc{\Tsupp}{\Supp^{\mathrm E}}

\nc{\abu}{\Vec{\bu}}
\nc{\abv}{\Vec{\bv}}
\nc{\av}{\Vec{v}}
\nc{\au}{\Vec{u}}

\nc{\FN}{\mathrm{FN}}
\nc{\DFN}{\mathrm{DFN}}

\nc{\ake}{\hs{.2ex}}

\nc{\va}{\Vec{a}}
\nc{\qtq}[1][and]{\quad\text{#1}\quad}
\nc{\qt}[1]{\quad\text{#1}}
\nc{\Mat}{\on{Mat}}
\nc{\ms}{\mspace}
\nc{\sprec}[1][\theta]{\mathop\prec\limits_{\,%
\raisebox{-.2ex}{$\scriptstyle\e^{i\,{#1}}$}}}
\nc{\spreceq}[1][\theta]{\mathop\preceq\limits_{
\raisebox{-.2ex}{$\scriptstyle\e^{i\,{#1}}$}}}
\nc{\Opc}[1][M]{\mathrm{Op}^{\mathrm{sub,\ms{2mu}c}}_{#1}}
\nc{\indd}{\inddlim}
\nc{\sindd}{\sinddlim}
\nc{\Ms}[1][M]{{#1}_{\sa}}

\newcommand{\bordered}[1]{{\mathsf{#1}}}
\newcommand{\bopen}[1]{{#1}}
\newcommand{\bclose}[1]{{\accentset{\vee}{#1}}}
\newcommand{\unbordered}[1]{{\accentset{\circ}{#1}}}

\nc{\unb}{\unbordered}
\newcommand{\oM}{\bopen{M}}
\newcommand{\cM}{\bclose{M}}

\newcommand{\oN}{\bopen{N}}
\newcommand{\cN}{\bclose{N}}
\newcommand{\bM}{\bordered{M}}
\newcommand{\bN}{\bordered{N}}

\nc{\bX}{{\widehat X}}
\nc{\cX}{{X^{\mathrm{c}}}}
\nc{\bY}{\widehat Y}
\nc{\bL}{\widehat L}
\nc{\bR}{{\ol\R}}
\nc{\bV}{{\ol \VV}}
\nc{\bW}{{\ol \WW}}

\nc{\olG}[1][f]{{\overset{\ms{4mu}\rule[-.05ex]{1.6ex}{.115ex}}{\Gamma}}_{%
\ms{-3mu}#1}}
\nc{\cf}{\bclose{f}}
\newcommand{\unbM}{\unbordered{\bM}}
\newcommand{\unbN}{\unbordered{\bN}}

\nc{\sub}{{\raisebox{.3ex}{$\scriptstyle\ms{3mu}\mathrm{sub}$}}}
\nc{\sC}{\C^\sub}
\nc{\sA}{\sha^\sub}
\nc{\sD}{\D^\sub}
\nc{\bs}{bordered subanalytic space\xspace}
\nc{\bss}{bordered subanalytic spaces\xspace}
\nc{\al}{\alpha}
\nc{\dbar}{{\ol\partial}}
\nc{\pDs}[1][X]{{(\pi_{#1}^{-1}\D_{#1})^\sub}}
\nc{\OT}[1][X]{{\sho_{#1}^\Tmp}}
\nc{\Monoto}[2]{\xymatrix{{#1}\ar@{ >->}[r]&{#2}}}
\nc{\OmT}[1][X]{\Omega^\Tmp_{#1}}
\nc{\epi}[1][M]{\pi_{#1}}
\nc{\piX}[1][X]{\pi_{#1}}
\nc{\tS}{{\mathscr{S}}^{\ms{2mu}\enh}}
\nc{\twS}{\widetilde{\mathscr{S}}^{\ms{5mu}\enh}}
\nc{\twD}{\widetilde{D}}
\nc{\bc}{\begin{cases}}
\nc{\ec}{\end{cases}}
\nc{\irc}[1][M]{\iota_{#1}^{\ms{4mu}\Rc}}

\begin{document}

\title
[{Riemann-Hilbert correspondence}]
{Riemann-Hilbert correspondence for
irregular holonomic D-modules}

\author[M.~Kashiwara]{Masaki Kashiwara}
\address{Research Institute for Mathematical Sciences\\
Kyoto University\\
Kyoto 606-8502, Japan}
\email{masaki@kurims.kyoto-u.ac.jp}

\thanks{The research
was supported in part by Grant-in-Aid for Scientific Research (B)
15H03608, Japan Society for the Promotion of Science.}

\date{}

\maketitle

\begin{abstract}
This is a survey paper on the Riemann-Hilbert correspondence
on (irregular) holonomic D-modules,
based on the 16-th Takagi lecture (2015/11/28).
In this paper, we use subanalytic sheaves,
an analogous notion to the one of indsheaves.
\end{abstract}
\tableofcontents

\section*{Introduction}
The classical Riemann-Hilbert problem asks for the existence of a linear  ordinary
differential equation with regular singularities and
a given monodromy on a curve.

Pierre Deligne (\cite{De70}) formulated it as a correspondence between 
integrable connections with regular singularities on a complex manifold $X$
with a pole on a hypersurface $Y$ and local systems
 on $X\setminus Y$.

Later the author constructed an equivalence of triangulated categories
between
$\Derb_\rh(\D_X)$, the derived category
of $\D_X$-modules 
with regular holonomic cohomologies,
and $\Derb_\Cc(\C_X)$, the derived category of sheaves on $X$ with 
$\C$-constructible cohomologies (\cite{Ka80, Ka84}).
The equivalence is given by the solution functor
$$\sol_X\cl \Derb_\rh(\D_X)\isoto \Derb_\Cc(\C_X)^\op.$$
Here $\sol_X(\shm)=\rhom[\D_X](\shm,\OO)$.
Note that $\Derb_\rh(\D_X)$ is self-dual by the duality functor.

\medskip
However, it was a long-standing problem to generalize it
to the (not necessarily regular) holonomic D-module case.
One of the difficulties
was that we could not find an appropriate substitute of the target category $\Derb_\Cc(\C_X)$.
Recently, the author solved it 
jointly with Andrea D'Agnolo (\cite{DK13})
by using an enhanced version of indsheaves.

There are two ingredients for the solution.

One is the notion of indsheaves.
This notion was introduced with Pierre Schapira in \cite{KS01}
to treat ``sheaves'' of functions with tempered growth,
such as $\Dbt$ of tempered distributions 
or $\Ot[]$ of tempered holomorphic functions.

The other ingredient is adding an extra variable. 
We consider indsheaves on $M\times \R$, not
on the base manifold $M$.
This method was originally introduced by Dmitry Tamarkin (\cite{Ta08}) 
in order to
treat non homogeneous Lagrangian submanifolds of the cotangent bundle 
in the framework of sheaf theory. In our context,
this method affords an appropriate language to capture various
growth of solutions at singular points.

Among the results used in the course of the proof 
is the description of the structure of flat connections due to 
Takuro Mochizuki (\cite{Mo09,Mo11}) and Kiran Kedlaya (\cite{Ke10,Ke11}).

 In this survey paper, we explain an outline of 
the irregular Riemann-Hilbert 
problem. We use here, instead of the notion of
indsheaves,
the analogous notion of ``subanalytic sheaves''.

For a complex manifold $X$, we construct a triangulated category
$\Derb(\sC_{X\times\fR})$, called the triangulated category of enhanced
subanalytic sheaves, a fully faithful functor
$e\cl \Derb(\C_X)\to \Derb(\sC_{X\times\fR})$ and its left quasi-inverse
$\fhom(\C_X^\enh,\scbul)\cl \Derb(\sC_{X\times\fR})\to \Derb(\C_X)$.
Next we construct $\OEn_X\in \Derb(\sC_{X\times\fR})$, 
the enhanced subanalytic sheaf of tempered holomorphic functions
such that $\fhom(\C_X^\enh,\OEn_X)\simeq\OO$.
By using $\OEn_X$ instead of $\OO$,
we define the enhanced solution functor
from the bounded derived category
$\BDC(\D_X)$ of $\D_X$-modules to
the category $\Derb(\sC_{X\times\fR})$ of enhanced subanalytic sheaves by
\eqn
\solE_X(\shm) \seteq \rhom[\D_X](\shm,\OEn_X)
\qt{for $\shm\in\BDC(\D_X)$.}
\eneqn

Restricting it to $\BDC_\hol(\D_X)$, the subcategory of
$\BDC(\D_X)$ consisting of complexes with holonomic cohomologies, 
we obtain a fully faithful functor
$$\solE_X\cl \BDC_\hol(\D_X)\monoto \Derb(\sC_{X\times\fR})^\op.$$
Furthermore, we have an isomorphism 
$$\fhom\bl\solE_X(\shm),\OEn_X\br\simeq\shm\qt{for any $\shm\in\BDC_\hol(\D_X)$.}
$$ 
Thus we obtain a quasi-commutative diagram:
\[
\xymatrix@C=4.5em{
\BDC_\reghol(\D_X)\ar[r]_{\Sol}^\sim \ar@{ >->}[d]\ar@/^2pc/[rrr]^\id
& \BDC_\Cc(\C_X)^\op \ar[rr]_{\rhom(\ast\,,\,\Ot)}^\sim \ar@{ >->}[d]^e
&& \BDC_\reghol(\D_X)\ar@{ >->}[d] \\
\BDC_\hol(\D_X) \ar[r]_{\solE_X} \ar@/_2.7pc/[rrr]^{\txt{\scriptsize{canonical}}}
& \Derb(\iC_{X\times\fR})^\op \ar[rr]_{\fhom(\ast\,,\,\OEn_X)}
&& \BDC(\D_X).
}
\]

\medskip
This paper is organized as follows.
In the first section, we review the local theory of
linear ordinary differential equations.
In the next sections, we shall review sheaves, D-modules and subanalytic sheaves.
After introducing the subanalytic sheaves of tempered distributions
and that of tempered holomorphic functions in \S\,\ref{sec:tempered},
we define the enhanced version of the de Rham functor and solution functor.
Then, in \S\,\ref{sec:Main}, we state our main theorems by using these functors.
In the next section \S\,\ref{sec:Out}, we give a very brief outline of the proof of the main theorems
by using the results of T.\ Mochizuki and K.\ S.\ Kedlaya.

In the last section \S\,\ref{sec:Stokes},
we explain how Proposition~\ref{prop:Stokes} 
on the Stokes phenomena in the one-dimensional case
can be interpreted in terms of the enhanced solution functors.

\vs{3ex}
We refer the reader to
\cite{DK13, KS14, KS15, DK15} for a more detailed theory.
Remark that the description of the Riemann-Hilbert correspondence
in this paper is different from that of
loc.\ cit.\ in the following points.
\bna
\item
We use in loc.\ cit.\ indsheaves instead of subanalytic sheaves.
Since the category of subanalytic sheaves can be embedded into that of
indsheaves, these two descriptions are almost equivalent.
\item In loc.\ cit., the category $\mathsf{E}^{\mathrm{b}}(\mathrm{I}\C_M)$ of 
enhanced indsheaves is defined as a quotient category 
of the category $\Derb(\mathrm{I}\C_{M\times\fR})$ of indsheaves on $M\times\fR$.
However, $\mathsf{E}^{\mathrm{b}}(\mathrm{I}\C_M)$
 can be also embedded into 
$\Derb(\mathrm{I}\C_{M\times\fR})$ by the right adjoint 
$\mathrm{R}^{\mathsf{E}}\cl \mathsf{E}^{\mathrm{b}}(\mathrm{I}\C_M)\to 
\Derb(\mathrm{I}\C_{M\times\fR})$
of the quotient functor. In our paper, we use the subanalytic sheaf version of
$\Derb(\mathrm{I}\C_{M\times\fR})$
instead of  $\mathsf{E}^{\mathrm{b}}(\mathrm{I}\C_M)$ by using the embedding 
$\mathrm{R}^{\mathsf{E}}$.
\ee

\section{Linear ordinary differential equations}
\label{sec:ordinary}
\subsection{One dimensional case}

Let us recall the local theory of linear ordinary differential equations.
Let $X\subset \C$ be an open subset with $0\in X$
and let $\shm$ be a holonomic $\shd_X$-module
such that $\SSupp(\shm)\subset \{0\}$ and
$\shm\simeq\shm(*\{0\})\seteq\OO(*\{0\})\tens[\OO]\shm$. 
Here $\OO(*\{0\})$ is the sheaf of meromorphic functions 
with possible poles at $0$.
It is equivalent to saying that
$\shm$ is a $\D_X$-module which is locally isomorphic to
$\OO(*\{0\})^{r}$ for some $r\in\Z_{\ge0}$ as an $\OO$-module.
Let us take a system of generators
$\{u_1,\ldots, u_r\}$ of $\shm$ as a free $\OO(*\{0\})$-module
on a neighborhood of $0$.
Then, writing $\abu$ for the column vector with these generators as components, 
we have
\eq
&&\dfrac{d}{dz}\abu=A(z)\abu
\label{eq:def}
\eneq
for some $A(z)\in\Mat_r\bl\sho_X(*\{0\})\br$, i.e.,
for an $(r\times r)$-matrix $A(z)$ whose components are in $\OO(*\{0\})$.
Then for any $\D_X$-module $\shl$ such that $\shl\simeq\shl(*\{0\})$,
we have
\eqn
&&\ba{ll}
\hom[\D_X](\shm,\shl)=\{\au\in\shl^r\;;\;&\text{$\au$ satisfies 
the same differential}\\
&\hs{1ex}\hfill \text{equation as \eqref{eq:def}\},}
\ea
\label{eq:solM}
\eneqn
where we associate  to $\au$ the morphism from $\shm$ to $\shl$
defined by $\abu\mapsto\au$.

\subsection{Regular singularities}
If we can choose a system of generators $\{u_1,\ldots, u_r\}$ of $\shm$
such that $zA(z)$ has no pole at $0$,
then we say that $0$ is a regular singularity of $\shm$,
or $\shm$ is  regular.
In such a case, 
there are  $r$ linearly independent solutions of the form
\eqn
\au_j=z^{\la_j}\sum_{s=0}^{r-1}\va_{j,s}(z)(\log z)^s\quad (j=1,\ldots,r),
\eneqn
where $\va_{j,s}(z)$ is a vector of holomorphic functions
defined on a neighborhood of $0$.
Hence,  after a change of generators
$\abv=D(z)\abu$ with some invertible matrix $D(z)\in\GL_r\bl\sho_X(*\{0\})\br$,
the new variable $\abv$ satisfies the equation
\eqn z\dfrac{d}{dz}\abv=C\,\abv
\eneqn
for some constant matrix $C\in\Mat_r(\C)$.
Then, by reducing $C$ to a Jordan form, we see that
$\shm$ is isomorphic to a direct sum of $\D_X$-modules
$\D_X(*\{0\})/\D_X(*\{0\})(z\dfrac{d}{dz}-\lambda)^{m+1}$ 
 with $\la\in\C$ and $m\in\Z_{\ge0}$.
Note that
$$\D_X(*\{0\})/\D_X(*\{0\})(z\dfrac{d}{dz}-\lambda)^{m+1}
\simeq\D_X/\D_X(z\dfrac{d}{dz}-\lambda-k)^{m+1}$$
for any $k\in\Z$ such that $\la+k\not\in\Z_{\ge0}$.

Recall that the solution sheaf of $\shm$ is defined by
$$\Sol(\shm)\seteq\rhom[\D_X](\shm,\OO).$$
Then the local system on $X\setminus\{0\}$
\eq 
&&L\seteq\Sol(\shm)\vert_{X\setminus\{0\}}
=\set{\au\in(\sho_{X\setminus\{0\}})^{\ms{0mu}r}}{\dfrac{d}{dz}\au=A(z)\au}
\eneq
has the monodromy $\exp(2\pi\sqrt{-1}C)$.
Hence $L$ completely determines $\shm$.

\subsection{Irregular singularities}
\label{subsec:Irr}
In the irregular case, we have the following results on the solutions of the
ordinary linear differential equation \eqref{eq:def}:
\bnum
\item
there exist linearly independent
$r$ formal  solutions $\hbu_j$ $(j=1,\ldots, r)$ of \eqref{eq:def}
with the form 
\eqn
\hbu_j=\ex^{\vphi_j(z)}z^{\la_j}\sum_{s=0}^{r-1}\va_{j,s}(z)(\log z)^s,
\eneqn
where $\vphi_j(z)\in z^{-1/m}\C[z^{-1/m}]$ for some $m\in\Z_{>0}$, $\la_j\in\C$, and
$$\text{$\va_{j,s}(z)=\sum\limits_{n\in m^{-1}\Z_{\ge0}}  \va_{j,s,n} z^n\in\C[[z^{1/m}]]^r$
with  $\va_{j,s,n}\in\C^r$,}$$
\item 
for any $\theta_0\in \R$ and each $j=1,\dots,r$, there exist an angular neighborhood
\eq
&&D_{\theta_0}=\set{z=r\ex^{i\,\theta}}{\text{$|\theta-\theta_0|<\vep$
and $0<r<\delta$}}\label{def:D}
\eneq
for sufficiently small $\eps,\delta>0$
and a holomorphic (column) solution $\bu_j\in\OO(D_{\theta_0})^r$
of \eqref{eq:def} defined on $D_{\theta_0}$
such that
$$\bu_j\sim \hbu_j,$$
in the following sense: for any $N>0$, there exists $C>0$ such that
\eq
&&|\bu_j(z)- \hbu_j^{N}(z)|\le C|\ex^{\vphi_j(z)}z^{\la_j+N}|
= C\ex^{\Re(\vphi_j(z))}|z^{\la_j+N}|,\label{eq:asymptotic}
\eneq
where $\hbu_j^{N}(z)$ is the finite partial sum
$$\hbu_j^{N}(z)=
\ex^{\vphi_j(z)}z^{\la_j}\sum_{s=0}^{r-1}
\hs{1ex}\sum\limits_{\substack{n\in m^{-1}\Z_{\ge0},\\[.2ex]\hs{-2ex}n\le N}}\hs{1ex} \va_{j,s,n}  z^n(\log z)^s.$$
Here we choose branches of $z^{1/m}$ and $\log z$ on $D_{\theta_0}$.
\ee

Note that a holomorphic solution $\bu_j$ is not uniquely determined
by the formal solution $\hbu_j$.
Indeed, $\bu_j+\sum_{k\not=j}c_k\bu_k$ also satisfies
the same estimate \eqref{eq:asymptotic}
whenever
\eqn
&&\text{$\Re(\vphi_k(z))<\Re(\vphi_j(z))$ on $D_{\theta_0}$
if $c_k\not=0$.}\label{eq:tri}
\eneqn

\subsection{Stokes phenomena}
We choose another sufficiently small angular domain $D_{\theta_1}$ 
such that $D_{\theta_0}\cap D_{\theta_1}\not=\emptyset$,
and, for each $j$, we take a holomorphic solution $\bu'_j$ defined on $D_{\theta_1}$
and with the asymptotic behavior \eqref{eq:asymptotic} on $D_{\theta_1}$.
Then we can write

$$\bu'_j=\sum_{k}a_{j,k}\bu_{k}\qt{on $D_{\theta_0}\cap D_{\theta_1}$}$$
with $a_{j,k}\in\C$.
Note that
\eq&&\text{
$\Re(\vphi_k(z))\le\Re(\vphi_j(z))$ on $D_{\theta_0}\cap D_{\theta_1}$
if $a_{j,k}\not=0$.}\label{cond:Stok}
\eneq

The matrix $(a_{j,k})_{1\le j,k\le r}$ is called the Stokes matrix.
If we cover a neighborhood of $\{0\}$ by such angular domains,
then a pair of adjacent angular domains gives a Stokes matrix,
and thus we obtain a family of matrices satisfying
\eqref{cond:Stok}.

Conversely,  we can find
a holonomic D-module $\shm$ whose Stokes matrices are
a given family of matrices  satisfying
\eqref{cond:Stok}.

\subsection{Stokes filtrations}\label{subsec:stoke}
\medskip
Deligne \cite{DMR07}
interpreted these results as follows
(see also Malgrange \cite{DMR07} and Sabbah \cite{Sa00,Sa13}).

Let $\vpi\cl\twX\to X$ be the real blow up of $X$ along $\{0\}$
defined in \S\;\ref{subsection:realblowup} below.
Namely,
\eqn
&\twX\seteq\set{(r,\zeta)\in\R_{\ge0}\times\C}{|\zeta|=1,
r\zeta\in X}\qtq\vpi(r,\zeta)=r\zeta.
\eneqn
Recall that
$L=\bl\Sol\shm\br\vert_{X\setminus\{0\}}$.
Let $S\seteq\vpi^{-1}(0)$ and $j\cl X\setminus\{0\}\to \twX$
and set
\eq
&&\twL=(\oim{j}L)\vert_S.
\eneq
Then $\twL$ is a local system on $S$ of rank $r$.

For the sake of simplicity, we assume that $m$ in  \S\,\ref{subsec:Irr} (i)
is equal to $1$.

Set $\Phi=\bl\sho_X(*\{0\})/\sho_X\br_0$.
For $\e^{i\,\theta_0}\in S$ and $\vphi,\psi\in \Phi$, we write
$\vphi\spreceq[{\theta_0}]\psi$ if there exists $c\in\R$ such that 
$\Re\bl\vphi(r\e^{i\,\theta})\br
\le\Re\bl\psi(r\e^{i\,\theta})\br+c$ for $0<r\ll1$ and $|\theta-\theta_0|\ll1$.
Then $\spreceq[\theta_0]$ is an order on $\Phi$.

For $\vphi\in \Phi$ and $\e^{i\,\theta}\in S$, we set
\eqn
&(F_\vphi)_{\e^{i\,\theta}}=\{u(z)\in(\twL)_{\e^{i\,\theta}}\;;\;
\parbox[t]{41ex}{$|u(z)|\le C\vert z^{-M}\e^{\vphi(z)}\vert$ on a neighborhood of $\e^{i\,\theta}$
for some $C>0$ and $M\in\Z_{>0}$\}.}
\eneqn
Then $\{F_\vphi\}_{\vphi\in\Phi}$ satisfies the following conditions
by the properties of the solutions explained in \S\,\ref{subsec:Irr}:
\bnum
\item $\{F_\vphi\}_{\vphi\in \Phi}$ is a filtration of $\twL$, namely,
\bna \item
$ F_\vphi$ is a subsheaf of $\twL$ for any $\vphi\in \Phi$, 
\item $\twL=\sum_{\vphi\in \Phi}F_\vphi$,
\item
$(F_{\vphi})_{\e^{i\,\theta}}\subset(F_\psi)_{\e^{i\,\theta}}$ if $\vphi\spreceq[\theta]\psi$,
\ee
\item for any $\e^{i\,\theta_0}\in S$, there exist an open neighborhood
$U$ of $\e^{i\,\theta_0}$, a finite subset $I$ of $\Phi$
 and a constant subsheaf $H_\vphi$ ($\vphi\in I$)
of $\twL\vert_U$
such that
\bna 
\item $\twL\vert_U=\soplus_{\vphi\in I}H_\vphi$,
\item for any $\e^{i\,\theta}\in U$ and $\vphi\in \Phi$, we have
$$(F_\vphi)_{\e^{i\,\theta}}=\soplus_{\substack{
\psi\in I,\;\psi\spreceq[{\theta}]\vphi}}(H_{\psi})_{\e^{i\,\theta}}.$$
\ee
\ee
If the above conditions are satisfied we say that
 $\{F_\vphi\}_{\vphi\in \Phi}$ is a Stokes filtration of the local system $\twL$.
Also in case $m>1$, we can define the notion of Stokes filtration
with a suitable modification.

\Prop\label{prop:Stokes}
The category of holonomic $\D_X$-module $\shm$
such that 
$$\text{$\SSupp(\shm)\subset \{0\}$ and
$\shm\simeq\shm(*\{0\})$ 
\ }$$
is equivalent to the category of pairs
$(L, \{F_\vphi\})$ of a local system $L$
on $X\setminus\{0\}$ and a Stokes filtration  $\{F_\vphi\}$ on
$\twL\seteq(\oim{j}L)\vert_S$.
\enprop
In order to generalize this result to holonomic D-modules in 
the several dimension case, we use enhanced subanalytic sheaves.
In the next sections, we shall review sheaves, D-modules and subanalytic sheaves.

\section{A brief review on sheaves and D-modules}\label{section:review}

\subsection{Sheaves}\label{subsec:shv}
We refer to~\cite{KS90} for all notions of sheaf theory used here.
For simplicity, we take the complex number field $\cor$ as the base field, 
although most of the results would remain true when $\cor$ is 
replaced with a commutative ring of finite global dimension.  

A topological space is {\em good} 
if it is Hausdorff, locally compact, countable at infinity and has finite flabby dimension.
 
One denotes by $\md[\cor_M]$  
the abelian category of sheaves of $\cor$-vector spaces 
on a good topological space $M$ and by 
$\Derb(\cor_M)$ its bounded derived category.
Note that $\md[\cor_M]$ has a finite homological dimension. 

For a locally closed subset $A$ of $M$, one denotes by $\cor_A$ the constant sheaf on $A$ with stalk $\cor$ extended by $0$ on $X\setminus A$. 

One denotes by $\Supp(F)$ the support of $F$.

There are many formulas concerning the six operations. 
For example, we have the formulas 
below in which $F,F_1,F_2\in\Derb(\cor_M)$, $G,G_1,G_2\in \Derb(\cor_N)$:
\eq
&&\ba{l}
\rhom(F_1\tens F_2,F)\simeq\rhom\bl F_1,\rhom(F_2,F)\br,\\
\roim{f}\rhom(\opb{f}G,F)\simeq\rhom(G, \roim{f}F),\\
\reim{f}(F\tens\opb{f}G)\simeq (\reim{f}F)\tens G\quad\text{
(projection formula),}\\
\epb{f}\rhom(G_1,G_2)\simeq\rhom(\opb{f}G_1,\epb{f}G_2),
\ea\label{eq:form1}
\eneq
and for a Cartesian square of good topological spaces,
\eqn\label{eq:CartSq}
\ba{c}\xymatrix@C=8ex{
M' \ar[r]^{f'} \ar[d]^{g'} & N' \ar[d]^{g} \\
M \ar[r]^{f}\ar@{}[ur]|-\square & N
}\ea
\eneqn
we have the {\em base change formulas}
\eq
&& \opb{g}\reim{f}\simeq \RR f'_{\ms{5mu}!}\ms{3mu}\opb{g'}\qtq\epb{g}\roim{f}\simeq 
\RR f'_{\ms{5mu}*}\ms{4mu}{g'}^{\ms{3mu}!}.
\label{eq:form2}
\eneq

\subsection{D-modules}\label{subsec:Dmod}
References for 
D-module theory are made to~\cite{Ka03}. See also \cite{Ka70,Ka75,Ka78,KK81,Bj93,HTT08}. 
Here, we shall briefly recall some basic constructions in  the theory of D-modules. 

Let $(X,\OO)$ be a {\em complex} manifold. 
We denote by 
\begin{itemize}
\item
$d_X$ the complex dimension of $X$,
\item
$\Omega_X$ the invertible $\OO$-module of differential forms of top degree,
\item
$\Omega_{X/Y}$ the invertible $\OO$-module $\Omega_X\tens[{\opb{f}\OO[Y]}]
\opb{f}(\Omega_Y^{\otimes-1})$ for a morphism $f\cl X\to Y$ of complex manifolds,
\item
$\Theta_X$ the sheaf of holomorphic vector fields,
\item 
 $\D_X$ the sheaf of algebras  of finite-order differential operators.
\end{itemize}

Denote by $\md[\D_X]$ 
the abelian category of left $\D_X$-modules and by  $\md[\D^\rop_X]$ that of  right $\D_X$-modules. There is an equivalence 
 \eq
 &&\hs{3ex}\mop\colon \md[\D_X] \isoto \md[\D_X^\op],\quad \shm\mapsto \shm^\mop\eqdot\Omega_X\tens[\OO]\shm. 
\eneq
By this equivalence, it is enough to study left $\D_X$-modules. 

The ring $\D_X$ is coherent and one denotes by $\mdc[\D_X]$ 
the thick abelian subcategory of $\md[\D_X]$ consisting of coherent modules. 
 
To a coherent $\D_X$-module $\shm$ one associates its characteristic variety $\chv(\shm)$, 
a closed $\C^\times$-conic 
{\em co-isotropic} (one also says {\em involutive}) $\C$-analytic subset of the cotangent bundle $T^*X$. 
The involutivity property is 
a central theorem of the theory and is due to~\cite{SKK73}. A purely algebraic proof was obtained later in~\cite{Ga81}.

If $\chv(\shm)$ is Lagrangian, $\shm$ is called {\em  holonomic}.
 It is immediately checked that the full subcategory $\mdhol[\D_X]$ of $\mdc[\D_X]$ consisting of holonomic $\D$-modules is a thick abelian subcategory.
 
A $\D_X$-module $\shm$ is \emph{quasi-good} 
if, for any relatively compact open subset 
$U\subset X$, there is a filtrant family $\{\shf_i\}_i$ of
coherent $(\OO\vert_U)$-submodules of $\shm\vert_U$
such that $\shm\vert_U=\sum_i\shf_i$. 
Here, a family $\{\shf_i\}_i$ is filtrant if, for any $i,i'$, there exists $i''$ such that $\shf_i+\shf_{i'}\subset \shf_{i''}$. 

A $\D_X$-module $\shm$ is  \emph{good} 
if it is quasi-good and coherent. The subcategories of $\md[\D_X]$ consisting of quasi-good (resp.\ good) $\D_X$-modules are abelian and thick. Therefore, one has  the triangulated categories
\eqn
&&\hs{-.5ex}\bullet\hs{.5ex}\BDC_\coh(\D_X)=\set{\shm\in\Derb(\D_X)}%
{\text{$H^j(\shm)$ is coherent for all $j\in\Z$}},\\
&&\hs{-.5ex}\bullet\hs{.5ex}\BDC_\hol(\D_X)=\set{\shm\in\Derb(\D_X)}%
{\text{$H^j(\shm)$ is holonomic for all $j\in\Z$}},\\
&&\hs{-.5ex}\bullet\hs{.5ex}\BDC_\rh(\D_X)=\set{\shm\in\Derb(\D_X)}%
{\text{$H^j(\shm)$ is regular holonomic for all $j\in\Z$}},\\
&&\hs{-.5ex}\bullet\hs{.5ex}\BDC_\qgood(\D_X)=\set{\shm\in\Derb(\D_X)}{\text{$H^j(\shm)$ is quasi-good  for all $j\in\Z$}},\\
&&\hs{-.5ex}\bullet\hs{.5ex}\BDC_\good(\D_X)=\set{\shm\in\Derb(\D_X)}{\text{$H^j(\shm)$ is good  for all $j\in\Z$}}.
\eneqn
 One may also consider the unbounded derived categories $\Der(\D_X)$, 
$\Der[+](\D_X)$   and $\Der[-](\D_X)$ and their full triangulated subcategories consisting of objects with coherent, holonomic,
regular holonomic, quasi-good and good cohomologies. 

We have the functors
\eqn
\rhom[\shd_X](\scbul,\scbul)&\cl& \Derb(\D_X)^\rop\times \Derb(\D_X)\to \Derp(\C_X),\\
\scbul\ltens[\D_X]\scbul&\cl&\Derb(\D_X^\rop)\times \Derb(\D_X)\to \Derb(\C_X).
\eneqn
We also have the functor
\eqn
\scbul\Dtens\scbul&\cl&\Derm(\D_X)\times\Derm(\D_X) \to\Derm(\D_X)
\eneqn
constructed as follows. 
For  $\D_X$-modules $\shm$ and $\shn$, the tensor product
$\shm\tens[\OO]\shn$ is endowed with a structure of
$\D_X$-module by
$$v(s\tens t)=(vs)\tens t+s\tens(vt)\qt{for $v\in\Theta_X$, $s\in\shm$ and
 $t\in\shn$.}$$
The functor $\scbul\Dtens\scbul$ is its left derived functor.
One defines the 
duality functor for D-modules by setting
\eqn
\hs{2ex}\Ddual_X\shm= \rhom[\D_X](\shm,\D_X\tens[\OO]\Omega_X^{\otimes-1})[d_X]
&\in&\Derb(\D_X)\quad\\
&&\hs{4ex}\text{for $\shm\in\Derb(\D_X)$.}
\eneqn

Now, let $f\cl X\to Y$ be a morphism of complex manifolds. 
The {\em  transfer  bimodule} 
$\D_{X\to Y}$ is a $(\D_X,\opb{f}\D_Y)$-bimodule defined as follows. As an 
 $(\sho_X,\opb{f}\D_Y)$-bimodule, 
$\D_{X\to Y}=\sho_X\tens[\opb{f}\sho_Y] \opb{f}\D_Y$. 
The left $\D_X$-module structure of 
 $\D_{X\to Y}$ is given by
\eqn
&&v(a\tens P)=v(a)\tens P+\sum_i aa_i\tens w_iP,
\eneqn
where  $v\in\Theta_X$ and  $\sum_ia_i\tens w_i$ is its image in 
$\sho_X\tens[\opb{f}\sho_Y]\opb{f}\Theta_Y$.

One also uses the opposite transfer  bimodule $\D_{Y\from X}=\opb{f}\D_Y\tens[{\opb{f}\OO[Y]}]\Omega_{X/Y}$, an $(\opb{f}\D_Y,\D_X)$-bimodule. 

Note that for another morphism of complex manifolds $g\cl Y\to Z$, one has the natural isomorphisms
\eqn
&&\D_{X\to Y}\ltens[\opb{f}\D_Y]\opb{f}\D_{Y\to Z}\simeq\D_{X\to Z},\\
&&\opb{f}\D_{Z\from Y}\ltens[\opb{f}\D_Y]\D_{Y\from X}\simeq\D_{Z\from X}.
\eneqn
One can now define the external operations on  D-modules by setting:
\eqn
&&\Dopb f\shn\eqdot \D_{X\to Y}\ltens[\opb{f}\D_Y]\opb{f}\shn
\qt{for $\shn\in\Derb(\D_Y)$,}\\
&&\Deim f\shm\eqdot\reim{f}(\shm\ltens[\D_X]\D_{X\to Y})
\qt{for $\shm\in\Derb(\D_X^\rop)$,}
\eneqn
and one defines $\Doim{f}\shm$ by replacing $\reim{f}$ with $\roim{f}$ 
in the  above formula. By using the opposite transfer bimodule $\D_{Y\from X}$ one defines similarly the inverse image of a right $\D_Y$-module or 
the direct images of a left $\D_X$-module. 

One calls respectively $\Dopb{f}$, $\Doim{f}$ and $\Deim{f}$ the inverse image,  direct image and proper direct image functors in the category of D-modules. 

Note that
\eqn
&&\Dopb f\OO[Y]\simeq\OO,\quad \Dopb f\Omega_Y\simeq\Omega_X.
\eneqn
Also note that the property of being quasi-good is stable by inverse image and tensor product, as well as by direct image by maps proper on the support of the module. The property of being good is stable by duality. 

\vspace{1ex}

Let $f\cl X\to Y$ be a morphism of complex manifolds. One associates the maps
\eqn\label{diag:tgmor}
\xymatrix@C=8ex{
T^*X\ar[rd]_-{\pi_X}&X\times_Y\ar[d]
\ar[l]_-{f_d}\ar[r]^-{f_\pi}T^*Y
                                      & T^*Y\ar[d]^-{\pi_Y}\\
&X\ar[r]^-f&Y.
}\eneqn
One says that $f$ is 
non-characteristic for $\shn\in\Derb_\coh(\D_Y)$ if the map $f_d$  is proper (hence, finite) on $\opb{f_\pi}\bl\chv(\shn)\br$.

The classical de Rham and 
solution functors are defined by
\begin{align*}
\dr_X &\cl \Derb(\D_X) \to \Derb(\C_X), &\shm &\mapsto \Omega_X \ltens[\D_X] \shm, \\
\sol_X &\colon \Derb(\shd_X)\to \Derb(\C_X)^\op, &\shm &\mapsto \rhom[\D_X] (\shm,\OO).
\end{align*}
For $\shm\in\Derb_\coh(\D_X)$, one has
\eq\label{eq:dualdrsolnt}
\sol_X(\shm) \simeq \dr_X(\Ddual_X\shm)[-d_X].
\eneq
Let us list up the relations of the de Rham functors with the 
inverse and direct image functors.

\begin{theorem}[{\rm Projection formulas~\cite[Theorems~4.2.8,~4.40]{Ka03}}]\label{th:Dprojform}
Let  $f\cl X\to Y$ be a morphism of complex manifolds. 
For $\shm\in\Derb(\D_X)$ and $\shl\in\Derb(\D_Y^\rop)$,
there are natural isomorphisms{\rm:}
\eqn
\Deim{f}(\Dopb{f}\shl\Dtens\shm)&\simeq&\shl\Dtens\Deim{f}\shm,\\
\reim{f}(\Dopb{f}\shl\ltens[\shd_X]\shm)&\simeq&\shl\ltens[\shd_Y]\Deim{f}\shm.
\eneqn
In particular, there is an isomorphism \lp commutation of the de Rham functor and direct images\rp
\eqn
\reim{f}(\dr_X(\shm))\simeq\dr_Y(\Deim{f}\shm).
\eneqn
\end{theorem}

\begin{theorem}[{\rm Commutation with duality~\cite{Ka03,Sc86}}]\label{th:oimopbDdual}
Let  $f\cl X\to Y$ be a morphism of complex manifolds.
\bnum
\item
Let  $\shm\in\Derb_\good(\D_X)$ and 
assume that $\Supp(\shm)$ is proper over $Y$. Then $\Deim{f}\shm\in\Derb_\good(\D_Y)$, and 
$\Ddual_Y(\Deim{f}\shm)\simeq\Deim{f}\Ddual_X\shm$.
\item
If $f$ is non-characteristic for 
$\shn\in\Derb_\coh(\D_Y)$, then 
 $\Dopb{f}\shn\in\Derb_\coh(\D_X)$ and
$\Ddual_X(\Dopb{f}\shn)\simeq\Dopb{f}\Ddual_Y\shn$.
\ee
\end{theorem}
\begin{corollary}\label{cor:Dadj2}
Let  $f\cl X\to Y$ be a morphism of complex manifolds. 
\bnum
\item
Let  $\shm\in\Derb_\good(\D_X)$ and  assume that $\Supp(\shm)$ is proper over $Y$. Then
we have the isomorphism for $\shn\in \Der(\D_Y)${\rm:}
\eqn
&&\hs{6ex}\roim{f}\rhom[\shd_X](\shm,\Dopb{f}\shn)\,[d_X]\simeq\rhom[\shd_Y](\Doim{f}\shm,\shn)\,[d_Y].
\eneqn
In particular, with the same hypotheses, we have the  isomorphism \lp commutation of the Sol functor and direct images\rp
\eqn
&&\hs{5ex}\roim{f}\Sol(\shm)\,[d_X]\simeq\Sol[Y](\Doim{f}\shm)\,[d_Y].
\eneqn
\item
 Let $\shn\in \Derb_\coh(\D_Y)$ and assume that $f$ is non-characteristic for $\shn$. Then
we have the isomorphism for $\shm\in \Der(\D_X)${\rm:}
\eqn
&&\hs{5ex}\roim f\rhom[\D_X](\Dopb f\shn,\shm)[d_X] \simeq
\rhom[\D_Y](\shn,\Doim f\shm)[d_Y].
\eneqn

\ee
\end{corollary}

\section{Subanalytic sheaves}

\subsection{Subanalytic spaces}

Let $M$ be a real analytic manifold. On $M$ there is the family of subanalytic subsets due to Hironaka (\cite{Hi73}) and Gabrielov (\cite{G68})
(see~\cite{BM88,VD98} for an exposition). 
This family is the smallest family of subsets of $M$  which satisfies the following properties:
\be[(a)]
\item for any real analytic manifold $N$ and
any proper morphism $f\cl N\to M$, the image of $N$ is subanalytic,
\item the intersection of two subanalytic subsets is subanalytic,
\item the complement of a subanalytic subset is subanalytic,
\item the union of a locally finite family of subanalytic subsets is subanalytic.
\ee
This family is a nice family. For example, it is closed by taking the closure and interior;
any relatively compact subanalytic subset has finitely many connected components,
and each connected component is subanalytic;
any closed subanalytic subset is the proper image of a real analytic manifold
as in (a).

\medskip
For real analytic manifolds $M$, $N$ and a closed subanalytic subset $S$ 
of $M$, we say that
a map $f\cl S\to N$ is subanalytic if its graph is subanalytic in 
$M\times N$. 
One denotes by $\sha^\R_S$ the sheaf of $\R$-valued subanalytic continuous maps on $S$.
A {\em subanalytic space}
$(M,\sha^\R_M)$, or simply $M$ for short,  is an $\R$-ringed space locally isomorphic to $(S,\sha^\R_S)$ 
for a closed subanalytic subset $S$ of a real analytic manifold.
In this paper, we assume that a subanalytic space is good, i.e.,
it is Hausdorff, locally compact, countable at infinity 
with finite flabby dimension. 

A morphism of subanalytic spaces is a morphism of $\R$-ringed spaces. 
Then we obtain the category of subanalytic spaces. 

We can define the notion of subanalytic subsets of a subanalytic space.

A sheaf $F$ on a subanalytic space $M$ is $\R$-\emph{constructible}
if there exists a locally finite family of locally closed subanalytic subsets 
$M_j$ ($j\in J$) such that
$M=\bigcup_{j\in J} M_j$ and
the sheaf $F\vert_{M_j}$ is locally constant of finite rank for each $j\in J$.
We denote by $\mdrc[\C_M]$ the full subcategory of
$\Mod(\C_M)$ consisting of $\R$-constructible sheaves.
It is a subcategory stable by taking kernels, cokernels and extensions.

One defines the category 
$\Derb_\Rc(\cor_M)$ as the full subcategory of $\Derb(\cor_M)$ consisting of objects $F$ such that $H^i(F)$ is $\R$-constructible for all  $i\in\Z$.
It is a triangulated subcategory and equivalent to
$\Derb\bl\mdrc[\C_M]\br$.

\subsection{Subanalytic sheaves}
Subanalytic sheaves are sheaves on a certain Grothendieck topology associated
with subanalytic spaces. Here we shall introduce it directly 
without using the language of Grothendieck topology.

Let $M$ be a subanalytic space.
Let $\Op$ be the category of open subsets.
The morphisms are inclusions, that is,
$\Hom[{\Op}](U,V)=\pt$ or $\emptyset$
according to $U\subset V$ or not.
Let $\Opc$  be
the full subcategory of $\Op$ consisting 
of relatively  compact subanalytic open subsets.

Recall that a sheaf is a contravariant functor from $\Op$ to $\Mod(\C)$
satisfying a certain ``patching condition''.
By replacing $\Op$ with $\Opc$ and modifying the ``patching condition'',
we obtain the notion of subanalytic sheaves introduced in
\cite{KS01} (see also \cite{Pr08} for its more detailed study).

\Def A subanalytic presheaf $F$ is a contravariant functor from
$\Opc$ to $\Mod(\C)$. We say that a  subanalytic presheaf $F$
is a subanalytic sheaf if it satisfies:
\bnum
\item
$F(\emptyset) = 0$,
\item
For $U,V\in\Opc$, the sequence
\[
0 \To F(U\cup V) \To[r_1] F(U) \dsum F(V) \To[r_2] F(U\cap V)
\]
is exact. Here $r_1$ is given by the restriction maps and $r_2$ is given by the difference of the restriction maps $F(U) \to F(U\cap V)$ and
$F(V) \to F(U\cap V)$.
\ee
\edf
Denote by $\Mod(\sC_{M})$ the category of subanalytic sheaves.
Recall that $\Mod(\C_{M})$ denotes the category of sheaves on $M$.
Since a sheaf is a contravariant functor from $\Op$,
the inclusion functor $\Opc\to\Op$ induces
a fully faithful functor
$$\iota_M\cl\Mod(\C_{M})\to \Mod(\sC_{M}).$$
For example,
$$\Hom[{\Mod(\sC_{M})}](\iota_M\C_U,F)\simeq F(U)\qt{for any $U\in\Opc$.}$$
The functor $\iota_M$ does not commute with inductive limits.
We denote by $\indd$ the inductive limit in $\Mod(\sC_{M})$
in order to avoid confusion.

Note that
$$\bl\inddlim[i]F_i\br(U)\simeq\indlim[i]\bl F_i(U)\br$$
for any $U\in\Opc$ and a filtrant inductive system $\{F_i\}_{i}$ of
subanalytic sheaves.

The functor $\iota_M$ admits a left adjoint, denoted by $\alpha_M$. For 
$F\in \md[\sC_M]$, the sheaf $\alpha_M(F)$ is the sheaf 
given by
$$\Op\ni U\longmapsto \prolim[{V\in\Opc,\ms{3mu}V\ssubset U}]F(V).$$
The functor $\al_M$ has a left adjoint $\beta_M$.
For $F\in\Mod(\C_M)$, $\beta_MF$ is the subanalytic sheaf
associated with the subanalytic presheaf
$\Opc\ni U\to F(\ol U)$.
Hence we have two pairs of adjoint functors
$(\al_M,\iota_M)$ and $(\beta_M,\al_M)$:
\eqn
&& \xymatrix{
{\md[\C_M]}\ar@<-1.7ex>[rr]_-{\beta_M}\ar@<1.7ex>[rr]^-{\iota_M}
&&{\md[\sC_M]}\ar[ll]|-{\;\al_M\;}.
}
\eneqn
Both $\Mod(\C_{M})$ and $\Mod(\sC_{M})$ are abelian categories,
and $\alpha_M$ and $\beta_M$ are exact.
The functor $\iota_M$ is left exact but {\em not} right exact.
However, we have the following result.
\Prop
The restriction of $\iota_M${\rm:}
\eq\irc\cl \mdrc[\C_M]\to \Mod(\sC_M)\label{def:irc}
\eneq
is exact.
\enprop

In fact, we have a more precise relation of these two categories
(see \cite{KS01}).

\Prop Let $\mdrcc[\C_M]$ be the category of $\R$-constructible sheaves on $M$ with compact supports.
Then, $\Mod(\sC_M)$ is equivalent to $\mathrm{Ind}(\mdrcc[\C_M])$,
the category of ind-objects in $\mdrcc[\C_M]$.
\enprop

For ind-objects we refer to \cite{SGA4} or \cite{KS06}.
In particular, we have
$$\Hom[ \Mod(\sC_M)](\iota_MG,\inddlim[i\in I]F_i)
\simeq\indlim[i\in I]\;\Hom[ \Mod(\sC_M)](\iota_MG,F_i)$$
for any $G\in\mdrcc[\C_M]$
and a filtrant inductive system $\{F_i\}_{i\in I}$ of subanalytic sheaves.

By the functor $\irc$, we regard $\R$-constructible sheaves as
 subanalytic sheaves.

We can define the restriction functor
$$\Mod(\sC_U)\to\Mod(\sC_V)\qtq[{for open subsets $U$ and $V\subset U$.}]
$$
For $F\in Mod(\sC_U)$, we denote by $F\vert_V\in\Mod(\sC_V)$ the image of $F$
by the restriction functor.

Hence, $\Op\ni U\mapsto \Mod(\sC_U)$ is a prestack on the topological space $M$.
\Prop
The prestack $\Op\ni U\mapsto \Mod(\sC_U)$ is a stack.
\enprop
We denote by $\hom$ the hom functor as a stack, i.e.,
for subanalytic sheaves $F_1, F_2$ on $M$, we define
\eqn
&&\sect(U;\hom(F_1, F_2))=\Hom[{\Mod(\sC_U)}](F_1\vert_U, F_2\vert_U)
\eneqn
for any open subset $U$ of $M$.
It is a sheaf on $M$.

\subsection{Bordered spaces}\label{sse:bord}
A \emph{bordered space} $\bM = (\oM,\cM)$ is a pair
of a good topological space $\cM$
and an open subset $\oM$ of $\cM$.

\begin{notation}
\label{not:Gammaf}
Let $\bM = (\oM,\cM)$ and $\bN = (\oN,\cN)$ be bordered spaces.
For a continuous map $f\colon \oM\to\oN$,
denote by $\Gamma_f \subset \oM\times\oN$ its graph, and by
$\olG$ the closure of $\Gamma_f$ in $\cM\times\cN$.
Consider the projections
\[
\xymatrix{
\cM & \ar[l]_-{q_1} \cM \times \cN \ar[r]^-{q_2} & \cN.
}
\]
\end{notation}

Bordered spaces form a category as follows:
a morphism $f\colon \bM \to \bN$ is a continuous map
$f\colon \oM\to\oN$ such that $q_1|_{\olG}\cl\olG\to\cM$ is proper;
the composition of two morphisms is the composition of the underlying continuous maps.

\begin{remark}\bnum
\item Let $f\cl M\to N$ be a continuous map.
\bna
\item
If $f$ can be extended to a continuous map $\cf\cl\cM\to\cN$,
then $f$ is a morphism of bordered space from $\bM$ to $\bN$.
\item If $\cN$ is compact, then $f$ is 
a morphism of bordered space from $\bM$ to $\bN$.
\ee
\item
The forgetful functor from the category of bordered spaces to that of good topological spaces is given by
\[
\bM = (\oM,\cM) \longmapsto \unbM \seteq \oM.
\]
It has a fully faithful left adjoint $M\mapsto(M,M)$. 
By this functor, we regard good topological spaces 
as particular bordered spaces, and denote $(M,M)$ 
simply by $M$.

Be aware that $\bM = (\oM,\cM) \mapsto \cM$ is not a functor.
\item 
Note that $\bM \simeq (\oM,\overline M)$, where $\overline M$ is the closure of $\oM$ in $\cM$.
More generally, for a morphism of bordered spaces $f\cl\bM\to\bN$,
$\bM$ is isomorphic to
the bordered space $(\Gamma_f,\olG)$.
\item The category of bordered spaces has an initial object,
the empty set. It has also a final object, $\pt$,
the topological space consisting of one point.
It  also admits products:
$$(\oM,\cM)\times (\oN,\cN) \simeq(\oM\times\oN,\cM\times \cN).$$
\ee
\end{remark}

Let $\bM = (\oM,\cM)$ be a bordered space.
The morphisms of bordered spaces
\begin{equation}
\label{eq:joM}
\oM \To \bM \To[j_\bM] \cM
\end{equation}
are defined by the
continuous maps $\oM \To[\id] \oM \hookrightarrow \cM$.

\begin{definition}
We say that a morphism $f\colon \bM \to \bN$ is \emph{semi-proper}
if $q_2|_{\olG}\cl\olG\to\cN$ is proper.
We say that $f$ is \emph{proper} if moreover $\unb f\cl\unbM\to\unbN$
is proper.
\end{definition}

For example, $j_\bM$ is semi-proper.

The class of semi-proper (resp.\ proper) morphisms is closed under composition.

\Def\label{def:subset}
A subset $S$ of a bordered space $\bM=(\oM,\cM)$ is a subset of $M$.
We say that $S$ is open (resp.\ closed, locally closed) if it is so in $\oM$.
We say that $S$ is relatively compact if it is contained
in a compact subset of $\cM$.
\edf

As seen by the following obvious lemma,
the notion of relatively compact subsets only depends on $\bM$
(and not on $\cM$).

\Lemma\label{lem:proper}
Let $f\cl \bM\to\bN$ be a morphism of bordered spaces.
\bnum
\item
If $S$ is a relatively compact subset of $\bM$, then its image
$\unb f(S)\subset \unbN$ is a relatively compact subset of $\bN$.
\item Assume furthermore that $f$ is semi-proper.
If $S$ is a relatively compact subset of $\bN$, then its inverse image
${\unb f}{}^{-1}(S)\subset \unbM$ is a relatively compact subset of $\bM$.
\ee
\enlemma

\subsection{Subanalytic sheaves on bordered subanalytic spaces}

A bordered subanalytic space is a bordered space $\bM=(\oM,\cM)$
such that $\cM$ is a subanalytic space and $\oM$ is a subanalytic open subset
of $\cM$.
Then we can consider the category of bordered subanalytic spaces.
A morphism $\bM=(\oM,\cM)\to\bN=(\oN,\cN) $ of bordered subanalytic spaces
is a morphism $f$ of bordered spaces such that the graph $\Gamma_f$
is a subanalytic subset of $\cM\times\cN$.

Let $\bM=(\oM,\cM)$ be a bordered subanalytic space.
We denote by $\Opc[\bM]$ the full subcategory of $\Op[\oM]$ consisting
of open subsets of $\oM$
which are subanalytic and relatively compact in $\cM$.
A subanalytic sheaf on $\bM$ is defined as follows.
\Def A subanalytic presheaf $F$ on 
a \bs $\bM$ is a contravariant functor from
$\Opc[\bM]$ to $\Mod(\C)$.
We say that  a subanalytic presheaf $F$ is a subanalytic sheaf if it satisfies:
\bnum
\item
$F(\emptyset) = 0$,
\item
For $U,V\in\Opc[\bM]$, the sequence
\[
0 \To F(U\cup V) \To[r_1] F(U) \dsum F(V) \To[r_2] F(U\cap V)
\]
is exact. 
\ee
\edf
We denote by $\Mod(\sC_\bM)$ the category of subanalytic sheaves on $\bM$.
We have a canonical fully faithful functor
\eq
&&\iota_\bM\cl\Mod(\C_{\unbM})\to\Mod(\sC_{\bM}).
\eneq
Here $\Mod(\C_{\unbM})$ denotes the category of sheaves on the topological space
$\unbM$. The functor $\iota_\bM$ is left exact but not exact.

We say that a sheaf on $\unbM$ is an $\R$-constructible sheaf on $\bM$ if
it can be extended to an $\R$-constructible
 sheaf on $\cM$. Let us denote by $\Mod_\Rc(\C_\bM)$ the category of
$\R$-constructible sheaves on $\bM$.
Then the restriction of $\iota_\bM$
$$\irc[\bM]\cl\Mod_\Rc(\C_{\bM})\to\Mod(\sC_{\bM})$$
is exact.
By this functor, we regard $\R$-constructible sheaves on
$\bM$ as subanalytic sheaves on $\bM$.

\subsection{Functorial properties of subanalytic sheaves}
\subsubsection{Tensor product and Inner hom}
Let $\bM=(M,\cM)$ be a \bs.
The category $\Mod(\sC_\bM)$ has tensor product and inner hom:
\eqn
\scbul\tens\scbul&\cl&\Mod(\sC_{\bM})\times \Mod(\sC_{\bM})\To\Mod(\sC_{\bM})
\qtq\\
\ihom(\scbul,\scbul)&\cl& \Mod(\sC_{\bM})^\op\times \Mod(\sC_{\bM})\To\Mod(\sC_{\bM}).
\eneqn
For $F_1, F_2\in\Mod(\sC_\bM)$, their tensor product
$F_1\tens F_2$ is the subanalytic sheaf associated with the subanalytic presheaf
$\Opc[\bM]\ni U\mapsto F_1(U)\tens F_2(U)$.
The inner hom $\ihom(F_1, F_2)$ is given by
$$\Opc[\bM]\ni U\mapsto \Hom[\Mod(\sC_{(U,\cM)})](F_1\vert_{(U,\cM)}, F_2\vert_{(U,\cM)}).$$
We have
$$\Hom[\Mod(\sC_\bM)](F_1\tens F_2, F_3)
\simeq\Hom[\Mod(\sC_\bM)](F_1,\ihom(F_2, F_3))
$$
for $F_1, F_2, F_3\in \Mod(\sC_{\bM})$.

The bifunctor $\scbul\tens\scbul$ is exact, and
$\ihom(\scbul,\scbul)$ is left exact.

\subsubsection{Direct images and Inverse images}
Let $\bM=(\oM,\cM)$ and $\bN=(\oN,\cN)$ be bordered subanalytic spaces
and let $f\cl \bM\to\bN$ be a morphism of bordered subanalytic spaces.

For $F\in \Mod(\sC_\bM)$, its direct image $\oim{f}F\in \Mod(\sC_\bN)$ is defined by
\eq
&&\hs{7ex}\bl\oim{f}F\br(V)=\Hom[\Mod(\sC_\bM)](\C_{f^{-1}V},F)
\qtq[{for any $V\in\Opc[\bN]$.}]
\eneq
The functor $\oim{f}\cl \Mod(\sC_\bM)\to\Mod(\sC_\bN) $ has a left adjoint
$$\opb{f}\cl \Mod(\sC_\bN) \to\Mod(\sC_\bM).$$
The functor $\opb{f}$ is called the inverse image functor.
For a subanalytic sheaf $G$ on $\bN$,
its inverse image $\opb{f}G$ is the subanalytic sheaf associated with the subanalytic presheaf
$$\Opc[\bM]\ni U\longmapsto\hs{-2ex} \sindlim[{V\in\Opc[\bN],\;U\subset f^{-1}V}]G(V).$$
The functor $\opb{f}$ is exact.

For $F\in \Mod(\sC_\bM)$, the direct image  with proper support
$\eeim{f}F$ is defined by
$$\sect(V; \eeim{f}F)=\indlim[{U}]\Hom(\C_{f^{-1}V};F\tens\C_U)
\qt{for $V\in\Opc[\bN]$.}$$
Here $U$ ranges over the open subsets in $\Opc[\bM]$ such that
$f^{-1}V\cap \ol U\to V$ is proper, where $\ol U$ denotes
the closure of $U$ in $M$.
In general, the diagram
$$\xymatrix@C=9ex
{\Mod(\C_\bM)\ar[r]^{\iota_\bM}\ar[d]^{\eim{f}}\ar@{}[dr]|{\mathrm{NC}}
&\Mod(\sC_\bM)\ar[d]^{\eeim{f}}\\
\Mod(\C_\bN)\ar[r]_{\iota_\bN}&\Mod(\sC_\bN)}
$$
is not commutative,
that is why we use the different notation $\eeim f$.
Note that the above diagram commutes if $f$ is semi-proper.

\Ex 
Let $M=\R_{>0}$, $N=\R$ and  let
$f\cl M\to N$ be the canonical inclusion. Then we have
\eqn
\eeim{f}\C_{M}&\simeq&\inddlim[{c\to0+}]\C_{\{t>c\}}\qtq\\
\eim{ f}\C_{M}&\simeq&\C_{\{t>0\}}.
\eneqn
They are not isomorphic.
Indeed, we have for $U=\set{t}{0<t<1}\in\Opc[N]$
$$\sect(U;\inddlim[{c\to\,0+}]\C_{\{t>c\}})\simeq\hs{-1ex}
\indlim[{c\to\,0+}]\sect(U;\C_{\{t>c\}})\simeq
0\qtq
\sect(U;\C_{\{t>0\}})\simeq\C.$$
\enex

Recall the morphism  $j_\bM\cl \bM\to\cM$ of \bss.
We have
\eq&&\ba{ll}
\eeim{j_\bM}\ms{1mu}\opb{j_\bM} F&\simeq\C_\unbM\tens F,\\
\oim{j_\bM}\ms{1mu}\opb{j_\bM} F&\simeq\ihom(\C_\unbM, F)
\ea\hs{3ex}\qt{for $F\in\Mod(\sC_{\cM})$.}
\eneq

Moreover, the functor $\opb{j_\bM}\cl\Mod(\sC_{\cM})\to \Mod(\sC_{\bM})$ induces an equivalence of abelian categories:
$$\Mod(\sC_{\cM})\;/\;\Mod(\sC_{\cM\setminus\oM})\simeq \Mod(\sC_{\bM}).$$
Here $\Mod(\sC_{\cM\setminus\oM})$ is regarded as a full subcategory of
$\Mod(\sC_{\cM})$ by the fully faithful exact functor $\oim{i}\simeq\eeim{i}
\cl \Mod(\sC_{\cM\setminus\oM})\to \Mod(\sC_{\cM})$,
where $i\cl\cM\setminus\oM\into \cM$ is the closed inclusion.

\subsection{Derived functors}

The fully faithful exact functor 
$$\irc[\bM]\cl\Mod_\Rc(\C_{\bM})\to\Mod(\sC_{\bM})$$
induces a fully faithful functor
$\Derb_\Rc(\C_{\bM})\monoto\Derb(\sC_{\bM})$
by which we regard $\Derb_\Rc(\C_{\bM})$ as a full subcategory of
$\Derb(\sC_{\bM})$.

The functors introduced in the previous subsection have derived functors:
\eqn
\scbul\tens\scbul&\cl&\Derb(\sC_{\bM})\times \Derb(\sC_{\bM})\To\Derb
(\sC_{\bM}),\\
\rihom(\scbul,\scbul)&\cl& \Derm(\sC_{\bM})^\op\times \Derp(\sC_{\bM})\To\Derp(\sC_{\bM}),\\
\opb{f}&\cl&\Derb(\sC_{\bN})\to \Derb(\sC_{\bM}),\\
\roim{f}&\cl&\Derb(\sC_{\bM})\to \Derb(\sC_{\bN}),\\
\reeim{f}&\cl&\Derb(\sC_{\bM})\to \Derb(\sC_{\bN}).
\eneqn
The functor $\reeim{f}$ has a right adjoint:
\eqn \epb{f}&\cl&\Derb(\sC_{\bN})\to \Derb(\sC_{\bM}).
\eneqn
If $\unb f\cl \unbM\to \unbN$ is topologically submersive,
i.e., $\unb f$ is isomorphic to $\unbN\times \R^n\to \unbN$
locally on $\unbM$, then
$$\epb f F\simeq\omega_{\unbM/\unbN}\tens\opb f F.$$
Here $\omega_{\unbM/\unbN}\seteq{\unb f}{}^{\ms{3mu}!}\ms{3mu}\C_{\unbN}\in\Derb_\Rc(\C_\bM)\subset \Derb(\sC_\bM)$
is the relative dualizing complex.

These six operations satisfy the  properties 
similar to \eqref{eq:form1} and \eqref{eq:form2}
for the Grothendieck's six operations for sheaves.

\subsection{Ring actions}\label{subsec:ring}
Let $M$ be a subanalytic space, and let $\sha$ be a sheaf of $\C$-algebras.
Let $F$ be a subanalytic sheaf on $M$.
We say that $F$ has an action of $\sha$,
or $F$ is a subanalytic $\sha$-module if a homomorphism of sheaves of
$\C$-algebras 
\eq\sha&&\to\hom(F,F)\label{eq:ringaction}\eneq is given.
Since $\hom(F,F)\simeq\al_M\ihom(F,F)$, the data \eqref{eq:ringaction}
is equivalent to
$$\beta_M\sha\to \ihom(F,F),$$
or $\beta_M\sha\tens F\to F$ with the associativity property.
 We denote by $\Mod(\sA)$ the category of subanalytic
$\sha$-modules, 
and by $\Derb(\sA)$ its bounded derived category.

We have the tensor functor and the hom functor:
\eqn
\scbul\tens[\sha]\scbul&\cl&\Derb(\sha^\op)\times\Derb(\sA)\to \Derm(\sC_M),\\
\rhom[\sha](\scbul,\scbul)&\cl&\Derb(\sha)^\op\times\Derb(\sA)\to \Derp(\sC_M).
\eneqn

\section{Subanalytic sheaves of tempered functions}
\label{sec:tempered}
\subsection{Tempered distributions}
Hereafter,
$M$ denotes a real analytic manifold.

An important property of subanalytic subsets is given by the lemma
below. (See Lojasiewicz~\cite{Lo59} and also~\cite{Ma66} 
for a detailed study of its consequences.)

\begin{lemma}\label{th:Loj1}
Let $U$ and $V$ be two relatively compact open subanalytic subsets of $\R^n$.
There exist a positive integer $N$ and $C>0$ such that
\eqn
&&\dist\big(x,\R^n\setminus (U\cup V)\big)^N
\le C \big(\dist(x,\R^n\setminus U)+\dist(x,\R^n\setminus V)\big).
\eneqn
\end{lemma}

We denote by $\Db_M$ 
the sheaf of Schwartz's distributions on $M$.
Denote by $\Dbt_M(U)$ the image 
of the restriction map $\sect(M;\Db_M)\to\sect(U;\Db_M)$,
and call it the space of {\em tempered distributions} on $U$. 

Using Lemma~\ref{th:Loj1}, one proves:
\Lemma
The subanalytic presheaf  $U\mapsto \Dbt_M(U)$ is a subanalytic sheaf on $M$. 
\enlemma
One denotes by $\Dbt_{M}$ this subanalytic sheaf. 
By the definition, there is a monomorphism
$$\Monoto{\Dbt_M}{\iota_M\Db_M},$$
and an isomorphism
$$\al_M\Dbt_M\simeq\Db_M.$$

Let us denote by $\D_M$ the sheaf of rings of differential operators with real analytic coefficients.
Then, $\Dbt_{M}$ is a subanalytic $\D_M$-module in the sense
of \S\;\ref{subsec:ring}.

\subsection{Tempered holomorphic functions}

Let $X$ be a complex manifold, and let us denote by $X_\R$ 
the underlying real analytic manifold.
We have defined the subanalytic sheaf of  tempered distributions
$\Dbt_{X_\R}$. It is a subanalytic $\D_{X_\R}$-module.
Let us consider the Dolbeault complex with coefficients in  $\Dbt_{X_\R}$:
$$\Dbt_{X_\R}\To[\dbar]\Omega_{\cX}^1\tens[{\OO[\cX]}]\Dbt_{X_\R}
\To[\dbar]\cdots \To[\dbar]\Omega_{\cX}^{d_X}\tens[{\OO[\cX]}]\Dbt_{X_\R}.
$$
Here $\cX$ is the complex conjugate manifold of $X$.
It is a complex in  the category
$\Mod(\sD_X)$ of subanalytic $\D_X$-modules.
Hence we can consider this complex as an object of $\Derb(\sD_X)$,
the bounded derived category of $\Mod(\sD_X)$.
We denote it by $\Ot$ and call it the subanalytic sheaf of tempered holomorphic functions.
Note that its cohomology groups are not concentrated at degree $0$ in general.

\subsection{Tempered de Rham and solution functors}
Setting $\Ovt_X\eqdot\Omega_X\tens[\sho_X]\Ot\in\Derb\bl(\D_X^\op)^\sub\br$,  we define 
the tempered de Rham and solution functors  by
\begin{align*}
\drt_X &\cl \Derb(\D_X) \to \Derm(\iC_X),&\shm &\mapsto \Ovt_X \ltens[\D_X] \shm, \\
\solt_X &\cl \Derb(\D_X)\to \Derp({\iC_X})^\op, &\shm& \mapsto \rhom[\D_X] (\shm,\Ot).
\end{align*}
One has
\[
\dr_X \simeq \alpha_X\circ\drt_X
\qtq
\sol_X \simeq \alpha_X\circ\solt_X.
\]
For $\shm\in\Derb_\coh(\D_X)$, one has
\eq\label{eq:dualdrsol}
\solt_X(\shm) \simeq \drt_X(\Ddual_X\shm)[-d_X].
\eneq

The next result is a reformulation of a theorem of~\cite{Ka84} (see also~\cite[Th.~7.4.1]{KS01})

\begin{theorem}\label{thm:ifunct0}
Let $f\cl X\to Y$ be a morphism of complex manifolds. 
There is an isomorphism in $\Derb\bl(\opb{f}\D^\rop_Y)^\sub\br$:
\eq\label{eq:funct1a}
&&\Ovt_X\ltens[\shd_X]\shd_{X\to Y}\,[d_X]\isoto \epb f \Ovt_Y\,[d_Y].
\eneq
\end{theorem}

  Note that this isomorphism \eqref{eq:funct1a} 
is equivalent to the isomorphism
\eqn
&&\hs{5ex}\shd_{Y\from X}\ltens[\shd_X]\Ot[X]\,[d_X]\isoto \epb f \Ot[Y]\,[d_Y]
\quad\text{in $\Derb\bl(\opb{f}\D_Y)^\sub\br$.}
\eneqn
\begin{corollary}\label{cor:ifunct1}
Let $f\cl X\to Y$ be a morphism of complex manifolds and let
$\shn\in\Derb(\D_Y)$. Then \eqref{eq:funct1a} induces  the isomorphism
\eqn
&&\drt_X(\Dopb f\shn)\, [d_X] \simeq \epb f \drt_Y(\shn) \,[d_Y]
\quad\text{in $\Derb(\iC_X)$.}\label{eq:funct1}
\eneqn
\end{corollary}

\Cor
For any complex manifold $X$, we have
$$\drt_X(\OO)\simeq\C_X[d_X].$$
\encor

The next results are a kind of Grauert direct image theorem
for tempered holomorphic functions,
and its  D-module version. 
\begin{theorem}[{{\rm Tempered Grauert theorem}\;{\cite[Th.~7.3]{KS96}}}]\label{th:tGrauert}
Let $f\cl X\to Y$ be a morphism of complex manifolds, let $\shf\in\Derb_\coh(\OO)$ and assume that $f$ is proper on $\Supp(\shf)$. Then there is a natural isomorphism
\eqn
&&\reeim{f}(\Ot\ltens[\OO]\shf)\simeq\Ot[Y]\ltens[{\OO[Y]}]\reim{f}\shf.
\eneqn
\end{theorem}

\Prop[{\cite[Th.~7.4.6]{KS01}}]\label{cor:ifunct2}
Let $f\cl X\to Y$ be a morphism of complex manifolds.
Let $\shm\in\Derb_\qgood(\D_X)$ and assume that $f$  is proper on $\Supp(\shm)$. 
Then there is an  isomorphism in $\Derb(\iC_Y)$
\eqn
\drt_Y(\Doim f\shm) &\isoto& \roim f\drt_X(\shm).\label{eq:funct2}
\eneqn
\enprop

For a  closed hypersurface $S\subset X$, denote by $\OO(*S)$
 the sheaf of meromorphic functions with poles at $S$. 
It is a holonomic $\D_X$-module and flat as an $\OO$-module. For $\shm\in\Derb(\D_X)$ or $\shm\in\Derb(\JD_X)$, set
\eqn
&&\shm(*S) = \shm \Dtens \OO(*S).
\eneqn

\Prop\label{le:GrotDR1}
Let $S$ be a closed complex hypersurface in $X$. There are isomorphisms
\eqn\ba{rcl}
\Ot(*S)&\simeq&\rihom(\C_{X\setminus S},\Ot)\qt{in $\Derb(\JD_X)$,}\\[1ex]
\OO(*S)&\simeq&\rhom
(\C_{X\setminus S},\Ot)\qt{in $\Derb(\D_X)$.}
\ea\eneqn
\enprop
\Cor\label{le:GrotDR2}
Let $S$ be a closed complex hypersurface in $X$. There are isomorphisms
in $\Derb(\sC_X)$
\eqn&&\ba{rl}
\drt_X(\OO(*S))&\simeq\dr_X(\OO(*S))\\
&\simeq\rhom(\C_{X\setminus S},\C_X)\,[d_X].
\ea
\eneqn
\encor

\medskip
\section{Enhanced subanalytic sheaves}
\subsection{Enhanced tensor product and inner hom}

Consider the 2-point compactification of the real line 
$\ol\R \eqdot \R
\sqcup\{+\infty,-\infty\}$. Denote by $\BBP^1(\R)=\R\sqcup\{\infty\}$ the real projective line. Then $\ol\R$ has a structure of subanalytic space such that the natural map $\ol\R\to\PR$ is a morphism of  subanalytic spaces.

\begin{notation}\label{not:Rinfty}
We will consider the \bs
\[
\R_\infty \eqdot (\R,\overline\R).
\]
\end{notation}

Note that $\R_\infty$ is isomorphic to $(\R,\PR)$ as a \bs.

Consider the morphisms of \bss
\begin{align}
a &\colon \R_\infty \to \R_\infty, \label{eq:muq1q2}\\
\mu, q_1,q_2 &\colon\R_\infty\times\R_\infty \to \R_\infty,\notag
\end{align}
where $a(t) = -t$, $\mu(t_1,t_2) = t_1+t_2$ 
and $q_1,q_2$ are the natural projections.

For a subanalytic space $M$, we will use the same notations for the associated morphisms
\begin{align*}
a &\colon M\times\R_\infty \to M\times\R_\infty, \\
\mu, q_1,q_2 &\colon M\times\R_\infty\times\R_\infty \to M\times\R_\infty.
\end{align*}
 We also use  the natural morphisms
\eq\ba{c}
\xymatrix@C=4ex{
M\times\R_\infty \ar[rr]\ar[dr]_\epi && M\times\overline\R \ar[dl]^{\overline\pi_M} \\
&M.
}\ea
\eneq

\begin{definition}\label{def:ctens1}
The functors
\begin{align*}
\ctens &\colon \BDC(\icor_{M\times\R_\infty}) \times \BDC(\icor_{M\times\R_\infty}) \to \BDC(\icor_{M\times\R_\infty}), \\
\cihom &\colon \Der[-](\icor_{M\times\R_\infty})^\op \times \Der[+](\icor_{M\times\R_\infty}) \to \Der[+](\icor_{M\times\R_\infty})
\end{align*}
are defined by
\begin{align*}
K_1\ctens K_2 &= \reeim {\mu} (\opb q_1 K_1 \tens \opb q_2 K_2), \\
\cihom(K_1,K_2)&= \roim {q_1} \rihom(\opb q_2 K_1, \epb{\mu}K_2).
\end{align*}
\end{definition}
One sets
\begin{align}\label{eq:tgeq0B}
\cor_{\{t\ge 0\}} &= \cor_{\{(x,t)\in M\times\R\;;\; t \ge 0\}}.
\end{align}
We use similar notation for  
$\cor_{\{t= 0\}}$, $\cor_{\{t> 0\}}$, $\cor_{\{t\le 0\}}$, $\cor_{\{t=a\}}$, etc.  
These are $\R$-constructible sheaves on $M\times\fR$.
We also regard them as objects of $\BDC(\icor_{M\times\R_\infty})$. 
\begin{lemma}
For $K\in\BDC(\icor_{M\times\R_\infty})$, there are isomorphisms
\[
\cor_{\{t= 0\}} \ctens K \simeq K \simeq \cihom(\cor_{\{t= 0\}}, K).
\]
More generally, for $a\in\R$, we have
\[
\cor_{\{t= a\}} \ctens K \simeq \roim{\mu_a} K \simeq \cihom(\cor_{\{t= -a\}}, K),
\]
where $\mu_a\colon M\times\R_\infty\to M\times\R_\infty$ is the morphism induced by the translation $t\mapsto t+a$.
\end{lemma}

\begin{corollary}
The category $\BDC(\icor_{M\times\R_\infty})$ has a structure of commutative tensor category with $\ctens$ as tensor product and $\cor_{\{t= 0\}}$ as unit object.
\end{corollary}

As seen in the following lemma, the functor $\cihom$ is the
inner hom of the tensor category $\BDC(\icor_{M\times\R_\infty})$.
\begin{lemma}
For $K_1,K_2,K_3\in\BDC(\icor_{M\times\R_\infty})$ one has
\eqn
&&\ba{ll}
&\Hom[\BDC(\icor_{M\times\R_\infty})](K_1\ctens K_2,\;K_3)\\[1ex]
&\hs{15ex}
\simeq \Hom[\BDC(\icor_{M\times\R_\infty})]\bl K_1,\;\cihom(K_2,K_3)\br,\\
&\cihom(K_1\ctens K_2,\;K_3)
\simeq \cihom\bl K_1,\;\cihom(K_2,K_3)\br,\\
&\roim\epi\rihom(K_1\ctens K_2,K_3) \simeq
\roim\epi\rihom(K_1,\cihom(K_2,K_3)).
\ea\label{eq:innerhom}
\eneqn
\end{lemma}

We define the outer hom functors on $\BDC(\icor_{M\times\R_\infty})$ as follows.
\begin{definition}\label{def:fihom}
One defines the hom functor 
\eqn&&
\fihom\colon \TDC(\icor_{M\times\fR})^\op \times \TDC(\icor_{M\times\fR}) 
\to \Derp(\icor_M)\\
&&\hs{20ex}\fihom(K_1,K_2) = \roim\epi\rihom(K_1,K_2),\
\eneqn
and one sets
\eqn
\fhom= \alpha_M\circ\fihom\,\cl\, \TDC(\icor_{M\times\fR})^\op \times \TDC(\icor_{M\times\fR}) \to \Derp(\cor_M).
\eneqn
\end{definition}
Note that
\begin{align*}
\Hom[\TDC(\icor_{M\times\fR})](K_1,K_2)
&\simeq H^0 \bl M;\fhom(K_1,K_2)\br.
\end{align*}

\subsection{Enhanced  sheaf of tempered distributions}

Let $M$ be a real analytic manifold.
Let $j_M\cl M\times\fR\to M\times\PR$ be the canonical morphism.

Let $t$ be the affine coordinate of $\PR$. Then,
$\partial_t\seteq\partial/\partial t$ is a vector field on $M\times\PR$, and hence it acts on $\Dbt_{M\times\PR}$.
\Lemma 
The morphism of subanalytic sheaves
$$\partial_t-1\cl \Dbt_{M\times\PR}\to\Dbt_{M\times\PR}$$
is an epimorphism.
\enlemma
We define the subanalytic sheaf on $M\times\fR$ by
$$\DbT_M=\ker\bl\partial_t-1:\opb{j_M}\Dbt_{M\times\PR}\to\opb{j_M}\Dbt_{M\times\PR}\br.$$ 
Since any solution of $(\partial_t-1)u(t,x)=0$ can be written as
$u(t,x)=\e^t\vphi(x)$,
we have a monomorphism in $\Mod(\sC_{M\times\fR})$
$$\Monoto{\DbT_M}{\opb{\epi}\iota_M\Db_M}\qt{by $u(t,x)\mapsto\vphi(x)$.}$$

Note that $\DbT_M$ is a subanalytic $\opb{\epi}\D_M$-module.
We call it the enhanced subanalytic sheaf of  tempered distributions.

\Prop
\eqn
\DbT_M&\simeq&\cihom(\C_{\{t\ge a\}}, \DbT_M)\qt{for any $a\in\R$}\\
&\simeq&\cihom(\C_M^\enh[1], \DbT_M).
\eneqn
Here we set
$$\C_M^\enh\seteq\inddlim[{c\to+\infty}]\C_{\{t<c\}}.$$
\enprop
The enhanced subanalytic sheaf $\C_M^\enh$ satisfies 
$$\C_M^\enh[1]\ctens \C_M^\enh[1]\simeq\C_M^\enh[1].$$
We can recover $\Dbt_M$ and $\Db_M$ from $\DbT_M$ as follows:
\eq\label{eq:oEot}
&&\ba{l}\fihom(\C_M^\enh,\DbT_M)\simeq\Dbt_M,\\
\fhom(\C_M^\enh,\DbT_M)\simeq\Db_M.
\ea
\eneq

\Rem The definition of $\DbT$ is slightly different from the one in 
\cite{DK13,KS15,DK15}.
The notation $\DbT$ in loc.\ cit.\ is equal to $\DbT[1]$ in our notation.
\enrem

\subsection{Enhanced  sheaf of tempered holomorphic functions}

Let $X$ be a complex manifold, and let us denote by $X_\R$ 
the underlying real analytic manifold.
We have defined the  enhanced subanalytic sheaf of  tempered distributions
$\DbT_{X_\R}$. It is a subanalytic $\opb{\pi_X}\D_{X_\R}$-module.
Let us consider the Dolbeault complex with coefficients in  $\DbT_{X_\R}$:
$$\DbT_{X_\R}\To[\dbar]\Omega_{\cX}^1\tens[{\OO[\cX]}]\DbT_{X_\R}
\To[\dbar]\cdots \To[\dbar]\Omega_{\cX}^{d_X}\tens[{\OO[\cX]}]\DbT_{X_\R}.
$$
Here $\cX$ is the complex conjugate manifold of $X$.
It is a complex in  the category
$\Mod(\pDs)$ of subanalytic $\opb{\pi_X}\D_X$-modules,
where $\DbT_{X_\R}$ is situated at degree $0$ and 
$\Omega_{\cX}^{d_X}\tens[{\OO[\cX]}]\DbT_{X_\R}$ at degree $d_X$.
Hence we can consider this complex as an object of $\Derb\bl\pDs\br$,
the bounded derived category of $\Mod\bl\pDs\br$.
We denote it by $\OT$ and call it the enhanced sheaf of tempered holomorphic functions.
Note that its cohomology groups are not concentrated at degree $0$.

\Rem\label{rem:exp}
If $X=\pt$, then
$$\OT\simeq\DbT_{X_\R}\simeq\C_X^\enh\seteq\inddlim[{c\to+\infty}]\C_{\{t<c\}}$$
as objects of $\Derb(\sC_\fR)$.
Indeed, for $-\infty\le a<b\le +\infty$,
$\e^t$ is a tempered distribution on the open interval $(a,b)$
if and only if $(a,b)\subset\{t<c\}$ for some $c\in\R$.

\enrem

By \eqref{eq:oEot}, we have
\eq
&&\ba{l}\fihom(\C_X^\enh,\OEn_X)\simeq\Ot\qtq\\
\fhom(\C_X^\enh,\OEn_X)\simeq\OO.
\ea\label{eq:oEOt}
\eneq

\subsection{Enhanced de Rham and solution functors}\label{subsection:Edrsol}
We set 
$$\OmT\seteq\opb{\pi_X}\Omega_X\tens_{\opb{\pi_X}\OO} \OT
\in\Derb\bl(\opb{\pi_X}\D_X^\op)^\sub\br.$$
We define the enhanced de Rham and solution functors
\begin{align*}
\drE_X &\colon \BDC(\D_X) \to \TDC(\iC_{X\times\fR}), \\
\solE_X &\colon \BDC(\D_X) \to \TDC(\iC_{X\times\fR})^\op
\end{align*}
by
\eqn
\drE_X(\shm) &\seteq& \OvE_X \ltens[\opb{\pi_X}\D_X]\opb{\pi_X} \shm, \\
\solE_X(\shm) &\seteq& \rhom[\opb{\pi_X}\D_X](\opb{\pi_X}\shm,\OEn_X).
\eneqn

Note that
\[
\solE_X(\shm) 
\simeq \drE_X(\Ddual_X\shm)[-d_X]\quad\text{for $\shm\in\Derb_\coh(\D_X)$.}
\]

By \eqref{eq:oEOt}, we have for any $\shm\in\Derb(\D_X)$
\eq
\ba{l}
\drt_X\shm\simeq \fihom(\C_X^\enh,\drE_X\shm),\\[1ex]
\dr_X\shm\simeq \fhom(\C_X^\enh,\drE_X\shm).
\ea\label{eq:DRo}
\eneq

For a particular case of holonomic D-modules, we can calculate
explicitly the enhanced de Rham.
Let $Y\subset X$ be a complex analytic hypersurface 
of a complex  manifold $X$,  and set $U=X\setminus Y$.
For $\varphi\in\OO(*Y)$, one sets
\begin{align*}
\D_X \ex^\varphi &= \D_X/\set{P\in\D_X}{P\ex^\varphi=0 \text{ on } U}, \\
\she^\varphi_{U|X}&=\D_X \ex^\varphi(*Y).
\end{align*}
Hence $\D_X \ex^\varphi$ is a $\D_X$-submodule of $\she^\varphi_{U|X}$, and $\D_X \ex^\varphi$ as well as
 $\she^\varphi_{U|X}$ is a holonomic $\D_X$-module.
Note that  $\she^\varphi_{U|X}$ is isomorphic to $\OO(*Y)$ as an $\OO$-module,
and the connection $\OO(*Y)\to\Omega^1_X\tens[\OO]\OO(*Y)$ is given by
$u\mapsto du+u d\vphi$. We call $\she^\varphi_{U|X}$ the exponential module
with exponent $\vphi$.

For $c\in\R$, write for short
\[
\{t<\Re \varphi+c\} \seteq \set{(x,t)\in U\times \R}{t<\Re \varphi(x)+ c} \subset X\times \R.
\]

Similarly to Remark~\ref{rem:exp}, one can calculate explicitly $\drE_X(\shm)$ when $\shm$ is an exponential D-module. 

\begin{proposition}\label{pro:Solphi2}
Let $Y\subset X$ be a closed complex analytic hypersurface, and set $U=X\setminus Y$.
For $\varphi\in\OO(*Y)$, there are isomorphisms 
\eqn
\drE_X(\she^\varphi_{U|X})&\simeq&\rihom(\opb{\piX}\C_U,\inddlim[c\to+\infty]\C_{\{t<\Re\varphi+c\}})[d_X].
\eneqn
\end{proposition}
The next results are easy consequences of Theorem~\ref{thm:ifunct0}, Corollary~\ref{cor:ifunct1},
Corollary~\ref{cor:ifunct2}.

\begin{theorem}\label{thm:Tfunct}
Let $f\colon X\to Y$ be a morphism of complex manifolds.
Let $f_\R\cl X\times\fR\to Y\times\fR$ be the morphism induced by $f$.
\bnum
\item
There is an isomorphism in $\TDC\bl(\opb{\pi_X}\opb f\D_Y)^\sub\br$
\[
\Tepb f \OEn_Y[d_Y] \simeq \opb{\pi_X}\D_{Y\from X} \ltens[\opb{\pi_X}\D_X] \OEn_X [d_X].
\]
\item
For any $\shn\in\BDC(\D_Y)$ there is an isomorphism in $\TDC(\iC_X)$
\eqn
\drE_X(\Dopb f \shn)[d_X] \simeq \Tepb f \drE_Y(\shn)[d_Y].
\eneqn
\item
Let $\shm\in\BDC_\good(\D_X)$, and assume that $\Supp(\shm)$ is proper over $Y$. 
Then, there are isomorphisms in $\TDC(\iC_Y)$ 
\eqn
&&\drE_Y(\Doim f\shm) \simeq  \Toim f\drE_X(\shm),
\eneqn
\ee
\end{theorem}

\section{Main theorems}\label{sec:Main}
The Riemann-Hilbert correspondence for holonomic D-modules can be stated as follows.
\Th\label{th:irrRH1}
There exists a canonical isomorphism functorial with respect to $\shm\in\Derb_\hol(\D_X)${\rm:}
\eq\label{eq:bjorkmorf1}
&&\shm\Dtens\Ot\isoto\fihom(\solE_X(\shm),\OEn_X)\text{\ in $\Derb(\iD_X)$.}
\eneq
\enth

Applying the functor $\alpha_X$ to~\eqref{eq:bjorkmorf1}, we obtain

\Th[{Enhanced Riemann-Hilbert correspondence}]\label{cor:irregRH2}
There exists a canonical isomorphism functorial with respect to 
$\shm\in\Derb_\hol(\D_X)${\rm:}
\eq\label{eq:bjorkmorf2}
&&\shm\isoto\fhom(\solE_X(\shm),\OEn_X)\ \text{in $\Derb(\D_X)$.}
\eneq
\enth

Thus we obtain the quasi-commutative diagram

\[
\xymatrix@C=4.5em{
\BDC_\reghol(\D_X)\ar[r]_{\Sol}^\sim \ar@{ >->}[d]\ar@/^2pc/[rrr]^\id
& \BDC_\Cc(\C_X)^\op \ar[rr]_{\rhom(\ast\,,\,\Ot)}^\sim \ar@{ >->}[d]^e
&& \BDC_\reghol(\D_X)\ar@{ >->}[d] \\
\BDC_\hol(\D_X) \ar[r]_{\solE_X} \ar@/_2pc/[rrr]_{\txt{\scriptsize{canonical}}}
& \Derb(\iC_{X\times\fR})^\op \ar[rr]_{\fhom(\ast\,,\,\OEn_X)}
&& \BDC(\D_X).
}
\]

Here the fully faithful functor
$e\cl \BDC(\C_X)\to\Derb(\iC_{X\times\fR})$ is defined by
$$e(F)\seteq\C^\enh_X\tens\opb{\piX}F.$$
Theorem~\ref{cor:irregRH2} shows that $\solE_X$ as well as $\drE_X$
is faithful.
In fact, we can also show the following full faithfulness of the enhanced de Rham
functor. 
\Th\label{th; emb}
For $\shm,\shn\in \Derb_\hol(\D_X)$,
one has an isomorphism
\eqn
&&\rhom[\D_X](\shm,\shn)\isoto \fhom(\drE_{X}\shm,\drE_{X}\shn).
\eneqn
In particular, the functor
\eqn
\drE_{X}:\Derb_\hol(\D_X)\To\TDC(\iC_{X\times\fR})
\eneqn
is fully faithful.
\enth

\begin{remark}Theorems~\ref{cor:irregRH2} and \ref{th; emb} 
due to~\cite[Th.~9.6.1, Th.~9.7.1]{DK13} are a natural formulation of 
the Riemann-Hilbert correspondence for irregular D-modules. 
Theorem~\ref{th:irrRH1} is due to~\cite[Th.~4.5]{KS14},
which is a generalization
of a theorem of J-E.~Bj\"ork (\cite{Bj93}). 
\end{remark}

\section{A brief outline of the proof of the main theorems}
\label{sec:Out}

We reduce the main theorems to the exponential D-module case,
using the results of Mochizuki and Kedlaya.

\subsection{Real blow up}\label{subsection:realblowup}
A classical tool in the study of differential equations is the
 real blow up.

Recall that $\C^\times$ denotes $\C\setminus\{0\}$ and $\R_{>0}$ the multiplicative group of positive real numbers. Consider the 
action of $\R_{>0}$ on $\C^\times\times\R$:
\eqn
&&\R_{>0}\times(\C^\times\times\R)\to\C^\times\times\R,\quad (a,(z,t))\mapsto (az,\opb{a}t)
\eneqn
and set
\eqn
&&\tw\C^\tot=(\C^\times\times\R)/ \R_{>0},\,
\tw\C=(\C^\times\times\R_{\ge0})/ \R_{>0},
\tw\C^{>0}=(\C^\times\times\R_{>0})/ \R_{>0}.
\eneqn
One denotes by  $\varpi^\tot$ the map:
\eq\label{eq:maptwom}
&&\varpi^\tot\cl \tw\C^\tot\to \C, \quad (z,t)\mapsto tz.
\eneq
Then we have
$$\tw\C^\tot\supset \tw\C\supset \tw\C^{>0}\isoto\C^\times.$$

Let $X=\C^n\simeq\C^r\times\C^{n-r}$ and let $D$ be the divisor $\{z_1\cdots z_r=0\}$, where $(z_1,\dots,z_n)$ is a coordinate system on $X$. Set 
\eqn
&&\twX^\tot=(\tw\C^\tot)^r\times\C^{n-r},\, 
\twX^{>0}=(\tw\C^{>0})^r\times\C^{n-r},\,
 \twX=(\tw\C)^r\times\C^{n-r}.
\eneqn
Then $\twX$ is the closure of $\twX^{>0}$ in $\twX^\tot$.
The map $\varpi^\tot$ in~\eqref{eq:maptwom} defines the map 
\eqn
&&\varpi\cl \twX\to X.
\eneqn
The map $\varpi$ is proper and induces an isomorphism
\eqn
&&\varpi\vert_{\twX^{>0}}\cl \twX^{>0}=\vpi^{-1}(X\setminus D)\isoto X\setminus D.
\eneqn

We call $\twX$ the {\em real blow up} of $X$ along $D$.

\begin{remark}
The real manifold $\twX$ (with boundary)
as well as the map $\varpi\cl \twX\to X$ may be intrinsically defined for a complex manifold $X$ and a normal crossing divisor $D$, but $\twX^\tot$ is only intrinsically defined as a germ of a manifold in a neighborhood of  $\tw X$. 
\end{remark}

\Def
Let
$\At$ be the subsheaf of
$\oim{j}(\OO[{X\setminus D}])$ consisting of holomorphic functions tempered
at any point of $\twX\setminus\twX^{>0}=\opb{\vpi}(D)$.
Here, $j\cl X\setminus D\simeq\twX^{>0}\into \twX$ is the open embedding.
We set 
\eqn&&\DA_\twX\eqdot\At\tens[\opb{\varpi}\OO]\opb{\varpi}\D_X.
\eneqn
\edf
 Then $\At$ and $\DA_\twX$ are sheaves of rings on $\twX$. We have a commutative diagram
$$\xymatrix@C=7ex@R=4ex{\opb\vpi\OO\ar@{ >->}[r]\ar@{ >->}[d]&\opb\vpi\D_X\ar@{ >->}[d]\\
\At\ar@{ >->}[r]&\DA_\twX.}
$$

We have
\eqn
&&\roim{\vpi}\At\simeq\OO(*D).\eneqn

For  $\shm\in \Derb(\D_X)$ we set:
\eq\label{eq:MA}
&&\shm^\tA\eqdot\DA_\twX\ltens[\opb{\varpi}\D_X]\opb{\varpi}\shm\in \Derb(\D^\At).
\eneq

Then we obtain
\eq
&&\roim{\vpi}\shm^\At\simeq\shm(*D).
\eneq

\subsection{Normal form}
The result in \S\;\ref{subsec:Irr} for ordinary linear differential equations
is generalized to higher dimensions
by T.\ Mochizuki (\cite{Mo09,Mo11}) and K.\ S.\ Kedlaya (\cite{Ke10,Ke11}).
In this subsection, we collect some of their results 
that we shall need. 

Let $X$ be a complex manifold and $D\subset X$ a normal crossing divisor. 
We shall use the notations introduced in the previous subsection:
in particular the real blow up $\vpi\cl\twX\to X$ and 
the notation $\shm^\tA$ of~\eqref{eq:MA}.

\begin{definition}\label{def:normal}
We say that a holonomic $\D_X$-module $\shm$ has \emph{a normal form} along $D$ if\\
${\rm (i)}$ $\shm\simeq \shm(*D)$,\\
${\rm (ii)}$  $\SSupp(\shm)\subset D$,\\
${\rm (iii)}$  for any $x\in\vpi^{-1}(D)\subset\twX$, there exist an open neighborhood $U\subset X$ of $\varpi(x)$ and finitely many $\varphi_i\in\sect(U;\sho_X(*D))$
such that
\eqn
&&(\shm^\tA)|_V \simeq
\left.\left( \bigoplus_i (\she_{U\setminus D|U}^{\varphi_i})^\tA \right)\right|_V
\eneqn
for some open neighborhood $V$ of $x$ with $V\subset\opb\varpi(U)$.
\end{definition}

A ramification 
of $X$ along $D$ on a neighborhood $U$ of $x\in D$ is a finite map
\[
p \colon X' \to U
\]
of the form 
$$p(z'_1,\ldots,z_n') = (z_1^{\prime\, m_1},\dots,z_r^{\prime\, m_r},z'_{r+1},\dots,z'_n)$$
for some $(m_1,\dots,m_r)\in(\Z_{>0})^r$.
Here $(z'_1,\dots,z'_n)$ is a local coordinate system of $X'$, and
$(z_1,\dots,z_n)$ is a local coordinate system of $X$
such that   $D=\{z_1\cdots z_r=0\}$.

\begin{definition}
\label{def:quasi-normal}
We say that a holonomic $\D_X$-module $\shm$ has \emph{a quasi-normal form} along $D$ if it satisfies (i) and (ii) in Definition~\ref{def:normal}, and if for any $x\in D$ there exists a ramification $p\colon X'\to U$ on a neighborhood $U$ of $x$ such that $\Dopb p (\shm|_U)$ has a normal form along $\opb p (D\cap U)$.
\end{definition}

\begin{remark}
In the above definition, 
$\Dopb p(\shm|_U)$ as well as $\Doim p\Dopb p(\shm|_U)$ is concentrated at degree zero. Moreover,
$\shm|_U$ is a direct summand of $\Doim p\Dopb p (\shm|_U)$.
\end{remark}

\subsection{Results of Mochizuki and Kedlaya}

The next result is an essential tool in the study of holonomic D-modules 
and is easily deduced from the 
fundamental work of  Mochizuki~\cite{Mo09,Mo11} 
(see also Sabbah~\cite{Sa00} for preliminary results  and see
Kedlaya~\cite{Ke10,Ke11} for the analytic case).

\begin{theorem}
\label{thm:normal}
Let $X$ be a complex manifold, $\shm$ a holonomic $\D_X$-module and $x\in X$.
Then there exist an open neighborhood $U$ of $x$, a closed analytic hypersurface $Y\subset U$, a complex manifold $X'$ and a projective morphism $f\colon X'\to U$ such that
\bnum
\item $\SSupp(\shm)\cap U\subset Y$,
\item $D\eqdot\opb f(Y)$ is a 
normal crossing divisor  of $X'$,
\item $f$ induces an isomorphism $X'\setminus D \to U \setminus Y$,
\item $(\Dopb f \shm)(*D)$ has a quasi-normal form along $D$.
\ee
\end{theorem}

Remark that, under assumption (iii), $(\Dopb f \shm)(*D)$ is concentrated at degree zero.

Using Theorem~\ref{thm:normal}, one easily deduces the next lemma.

\begin{lemma}\label{lem:redux}
Let $P_X(\shm)$ be a statement concerning a complex manifold $X$ and a holonomic object $\shm\in\BDC_\hol(\D_X)$. Consider the following conditions.
\bna
\item
Let $X=\Union\nolimits_{i\in I}U_i$ be an open covering. Then $P_X(\shm)$ is true if and only if $P_{U_i}(\shm|_{U_i})$ is true for any $i\in I$.
\item
If $P_X(\shm)$ is true, then $P_X(\shm[n])$ is true for any $n\in\Z$.
\item
Let $\shm'\to\shm\to\shm''\To[+1]$ be a distinguished triangle in $\BDC_\hol(\D_X)$. If $P_X(\shm')$ and $P_X(\shm'')$ are true, then $P_X(\shm)$ is true.
\item
Let $\shm$ and $\shm'$ be holonomic $\D_X$-modules. 
If $P_X(\shm\dsum\shm')$ is true, then $P_X(\shm)$ is true.
\item Let $f\colon X\to Y$ be a projective morphism and $\shm$ a good holonomic $\D_X$-module. If $P_X(\shm)$ is true, then $P_Y(\Doim f\shm)$ is true.
\item If $\shm$ is a holonomic $\D_X$-module with a normal form along a normal crossing divisor of $X$, then $P_X(\shm)$ is true.
\ee
If conditions {\rm (a)--(f)} are satisfied, 
then $P_X(\shm)$ is true for any complex manifold $X$ and any $\shm\in\BDC_\hol(\D_X)$.
\end{lemma}

\medskip
\noi
\emph{Sketch of the proof of the main theorems in \S\ake\ref{sec:Main}}\\
By applying Lemma \ref{lem:redux}, we reduce the assertions to the case of
holonomic D-modules with a normal form, then to the case of the  exponential D-modules.
\hfill\qedsymbol

\section{Stokes filtrations and enhanced de Rham functor}
\label{sec:Stokes}
In this last section, we explain the relation between 
the enhanced solution sheaf
and the Stokes filtration discussed in \S\;\ref{subsec:stoke}.
Let us keep the notations in \S\,\ref{subsec:Irr}.
In particular, recall that $0\in X\subset\C$, $\shm$ is
a holonomic $\D_X$-module, 
$\vpi\cl\twX\to X$ is the projection and
$j\cl X\setminus\{0\}\into \twX$ is the open embedding.
We set
$X^*\seteq X\setminus\{0\}$.
Let
$\vpi_\R\cl \twX\times\fR\to X\times\fR$ be the morphism induced by $\vpi$ 
and let $i\cl S\seteq\vpi^{-1}(0)\into \twX$ be the closed embedding.

Set $$\shm'=\Ddual_X\bl(\Ddual_X \shm)(*\{0\})\br.$$
Then we have a morphism $\shm'\to\shm$ such that
it induces an isomorphism $\shm'(*\{0\})\isoto\shm$.

We set 
\eqn&&
\tS\seteq\solE(\shm')
\simeq\rihom\bl\C_{X^*\times\R},\solE(\shm)\br
\in \Derb(\sC_{X\times\fR}).
\eneqn
Since $\e^{t}\bu_j(z)\vert_{D_{\theta_0}\times\R}$ is tempered on
$$\{t+\Re\vphi_j<c\}\seteq
\{(z,t)\in X^*\times\R\mathbin{;}t+\Re(\vphi_j(z))<c\}$$
for any $c$ (see Proposition~\ref{pro:Solphi2}),
we have
$$\C_{D_{\theta_0}\times\R}\tens\tS\simeq\bigoplus_{1\le j\le r}\C_{D_{\theta_0}\times\R}\tens
\tS_{\vphi_j},$$
where 
\eq
\tS_\vphi\seteq\inddlim[{c\to+\infty}]\C_{\{t+\Re\widetilde\vphi<c\}}\in
\Mod(\sC_{X\times\fR})\qt{for $\vphi\in\Phi$.}\label{eq:Sphi}
\eneq
Here $\widetilde\vphi\in(\sho_X(*\{0\})_0$ is a representative of
$\vphi\in \Phi\seteq(\sho_X(*\{0\}/\sho_X)_0$.
Note that the right-hand side of of \eqref{eq:Sphi} does not depend on the 
choice of a representaive $\widetilde\vphi$.

Set 
\eqn
\twS_\vphi&\seteq&\Topb{\vpi}\tS_\vphi
\in \Derb(\sC_{\twX\times\fR})
\qtq\\
\twS&\seteq&\Tepb{\vpi}\tS\simeq
\rihom\bl\C_{X^*\times\R},\Topb{\vpi}\tS\br\in \Derb(\sC_{\twX\times\fR}).
\eneqn
Then set 
$$K_\vphi\seteq\fhom(\twS_\vphi,\twS)\in\Derb(\C_\twX).$$
Since $\twS_\vphi\vert_{X^*\times \fR}\simeq\C_{X^*}^\enh$,
we have
$$K_\vphi\vert_{X^*}\simeq L\seteq\hom[\D_X](\shm,\OO)\vert_{X^*}.$$
Then we obtain a morphism of sheaves on $S$
$$\opb{i}K_\vphi\to \opb{i}\oim{j}(K_\vphi\vert_{X^*})\simeq\twL\seteq
i^{-1}j_*L.$$
\Lemma The object $\opb{i}K_\vphi\in\Derb(\C_S)$ is concentrated at degree $0$.
The above morphism 
$\opb{i}K_\vphi\to \twL$
is a monomorphism and its image coincides with $F_\vphi$.
\enlemma
Thus, $\solE(\shm)$ recovers the Stokes filtration $\{F_\vphi\}_{\vphi\in\Phi}$
on $\twL$.

\providecommand{\bysame}{\leavevmode\hbox to3em{\hrulefill}\thinspace}

\end{document}